\documentclass[a4paper,12pt,oneside]{amsart}
\usepackage{amssymb} \usepackage{xy} \xyoption{all} \CompileMatrices

\newtheorem{thm}{Theorem}[section] \newtheorem{prop}[thm]{Proposition}
\newtheorem{lemma}[thm]{Lemma} \newtheorem{cor}[thm]{Corollary}
 \newtheorem{dfn}[thm]{Definition}
\newtheorem{rem}[thm]{Remark}

\begin{document}

\title{The isometry group of Outer Space}

\author{Stefano Francaviglia} \email{francavi@dm.unibo.it}
\address{Dipartimento di Matematica Universit\`a di Bologna, P.zza
Porta S. Donato 5, 40126 Bologna (Italy) } \author{Armando Martino}
\email{A.Martino@soton.ac.uk} \address{School of Mathematics
University of Southampton University Road Southampton SO17 1BJ }

\begin{abstract} We prove analogues of Royden's Theorem for the Lipschitz
  metrics of Outer Space, namely that $Isom(CV_n)=Out(F_n)$.
\end{abstract}

\maketitle

\tableofcontents

\section{Introduction}\label{sintro} For $n\geq2$ let $F_n$ be the
free group of rank $n$, and $Out(F_n)$ be the group of outer
automorphisms of $F_n$. The Culler-Vogtmann {\em Outer Space}, $CV_n$,
is the analogue of Teichmuller space for $Out(F_n)$ and is a space of
metric graphs with fundamental group of rank $n$.

As for Teichmuller space, one can define the Lipschitz metric of
$CV_n$ with a resulting metric which is not symmetric. This
non-symmetric metric is geodesic and seems natural in terms of
capturing the dynamics of free group automorphisms; for instance the
axes of iwip automorphisms (\cite{Yael08}). However the non-symmetric
version also lacks some properties one might want; it fails to be
complete, for instance, while the symmetrised version turns $CV_n$
into a proper metric space (see~\cite{FrMa,Yael08,BY09}, and
also~\cite{U} for a different approach.)

The group $Out(F_n)$ naturally acts on CV$_n$
and the action is by isometries. It is also easy to see that this
action is faithful for $n \geq 3$ but not faithful for $n=2$. The
reason for this is that $Out(F_2) \simeq GL(2, \mathbb{Z})$ has a
central element of order $2$, namely $-I_2$, which is in the kernel of
the action. If one picks a basis, $x_1, x_2$ for $F_2$ the
automorphism which sends each $x_i$ to ${x_i}^{-1}$ is a pre-image in
$Aut(F_2)$ of $-I_2$.

In this paper, we prove an analogue of Royden's Theorem for both
metrics, and any
rank, so that $Isom(CV_n)=Out(F_n)$ (see below for exact statements).

There are many of this kind of results in literature, for instance
\begin{itemize}
\item The Fundamental Theorem of projective geometry (If a field $F$ has
  no non-trivial automorphisms, the group of
  incidence-preserving bijections of the projective space of dimension
  $n$ over F is precisely $PGL(n, F)$).
\item Tits Theorem: Under suitable hypotheses, the
full group of simplicial automorphisms of the spherical building
associated to an algebraic group is equal to the algebraic group
(\cite{Tits}).
\item Ivanov's Theorem: The group of simplicial
automorphisms of the curve-complex of a surface $S$ of genus at least two
is the mapping class group of $S$ (\cite{Ivanov}).
\item Royden's Theorem: The isometry group of the Teichmuller space of
  $S$ is the mapping class group of $S$ (\cite{Royden}).
\item Bridson and Vogtmann's Theorem: For $n\geq 3$ the group of
  simplicial automorphisms of the spine of $CV_n$ is $Out(F_n)$ (\cite{BrVo01}).
\item Aramayona and Souto's Theorem: For $n\geq 3$, the group of simplicial
automorphisms of the free splitting graph is $Out(F_n)$; (\cite{AS09}).
\end{itemize}

Our main results are:

\begin{thm}\label{tmain} With respect to the symmetric Lipschitz
distance,
$$Isom(CV_n)=Out(F_n) \ \mbox{\rm for } n\geq 3.$$

For $n=2$, $$Isom(CV_2)=PGL(2, \mathbb{Z}).$$
\end{thm}

We note that replacing the symmetric distance by its non-symmetrised
version one gets the same result.

\begin{thm}\label{tm2} For both non-symmetric Lipschitz distances
$d_R$ and $d_L$, Isom$(CV_n)$ is $Out(F_n)$ for $n \geq 3$ and $PGL(2,
\mathbb{Z})$ for $n=2$.
\end{thm}

This kind of result has immediate corollaries of fixed-point type
(see for example~\cite{BrFarb,BrVo01}).
\begin{cor}\label{c}
Let $G$ be a semisimple Lie group with finite centre and no compact factors
and suppose the real rank of $G$ is at least two.
Let $\Gamma$ be a non-uniform, irreducible lattice
in $G$. Then every isometric action of $\Gamma$ on $CV_n$ has a global
fixed point.
\end{cor}

\medskip

As we note above, there already exists a result of this kind for the
{\em spine} of $CV_n$, \cite{BrVo01}, which states that the
simplicial automorphism group of the spine of $CV_n$ is equal to
$Out(F_n)$ for $n \geq 3$. At a first glance, Theorem~\ref{tmain}
could appear to be a direct consequence of~\cite{BrVo01} after some easy
remarks (using, for instance, Lemma~\ref{L1}) and in fact that
was exactly the thought of the authors when this work started.

However, the main difficulty in the paper is precisely moving
from a statement that an isometry preserves the simplicial structure
of $CV_n$ to the statement that it is the identity. For instance, once
one knows that an isometry leaves some simplex invariant, it is not
clear, a priori, that the centre of the simplex is fixed (in fact it 
is not true in general if one simply looks at isometries of a simplex
rather than the restriction of a global isometry). And even when one
has that a given isometry leaves {\em every} simplex invariant, it is
not clear how to deduce that the isometry is in fact the identity -
obviously, this is in sharp contrast to the piecewise Euclidean
metric.

Let us emphasise this contrast. Suppose that one wants to prove
Theoorem~\ref{tmain} for the piecewise Euclidean metric. First, note
that simplices corresponding to graphs with disconnecting edges are an
obvious obstruction. However, one always may to restrict to a ``reduced''
Outer Space by removing such simplices. Now, looking at the incidence
structure of ideal vertices, one can prove that any isometry
(w.r.t. the piecewise Euclidean metric) maps ideal vertices to ideal
vertices and thus simplices to simplices because isometries are local
PL-maps. (Lemma~\ref{L1} is no longer
true, as stated, for this reduced Outer Space as one can easily see in
the rank-$2$ case.) Then, invoking the 
Bridson-Vogtmann result, one gets that up to composing with
automorphisms, simplices are not permuted, and the PL-structure now completes
the job. 

\medskip

Now let us return to the Lipschitz metric. There are four key facts in
Theorem~\ref{tmain}. First, the study of local 
isometries. The main point is that in general, the isometry group of a
fixed simplex of $CV_n$ is in fact much bigger than its stabiliser in
$Out(F_n)$.

The second fact is that $CV_n$ is highly
non-homogeneous. This allows one to find particular points in simplices of
$CV_n$ that
are invariant under isometries, so that one can characterise those
isometries that are restrictions of global ones.

Third, there is the fact that asymptotic behavior of distances from
a particular set of points determines the distance between points of
$CV_n$. This is a non-trivial issue, that we like to paraphrase
saying that {\em Busemann functions of ideal vertices are coordinates for
$CV_n$.} The main consequence of this fact is that one can
deduce that an isometry that does not permute
simplices is in fact the identity.

Lastly, there is a
permutation issue, similar to the one faced in~\cite{BrVo01}, that we
solve metrically using our ``Busemann functions''.

We also remark that Theorem~\ref{tmain} holds for any rank
and includes the study of simplices with disconnecting edges.
The complete schema of the proof of Theorem~\ref{tmain} is described
in Section~\ref{sschema}.

\medskip

\noindent{\bf Acknowledgments.} The first named author wishes to thank
the UFF of Niter\'oi (RJ, Brazil), the CRM of Barcelona (Spain)
and the SOTON (Southampton, UK) for
their kind hospitality and the great work environment they provide.

This work was inspired by the beautiful
articles~\cite{BrVo01,BrFarb}, and most of the material of this
introduction was picked from there.

\section{Preliminaries}

In this section we fix terminology, give basic definitions, and 
recall some known facts (and prove some easy ones) that we shall need
for the rest of the paper. Experienced readers may skip directly to 
next section and refer to present one just for notation.

\subsection{Outer Space} First of all, we recall what Culler-Vogtmann
space or ``Outer Space'' is. We refer to the pioneer
work~\cite{CuVo86} and beautiful surveys~\cite{Vog02,Vog08} for more
details.

For any $n\geq 2$ let $F_n$ be the free group of rank $n$ which we
identify with the fundamental group of $R_n=S^1\wedge\dots\wedge S^1$
(the product taken $n$ times).

Consider finite graphs $X$ whose vertices have valence at least three,
this means that each vertex has at least three germs of incident
edges.  We require that $X$ has rank $n$, that is to say,
$\pi_1(X)\simeq F_n$ and that $X$ comes equipped with a metric.
Giving a metric on $X$ is equivalent to giving positive lengths for
the edges of $X$.

We also require $X$ to be a {\em marked graph}, which is to say that
it comes with a fixed {\em marking}.  A marking on $X$ is a continuous
map $\tau:R_n\to X$ which induces an isomorphism $\tau_*:
F_n\simeq\pi_1(R_n)\to\pi_1(X)$. Two marked metric graphs $(A,\tau_A)$
and $(B,\tau_B)$ are considered equivalent if there exists a
homothety, $h:A \to B$, such that the following diagram commutes up to
free homotopy,

$$
\xymatrix{ A \ar[rr]^h & & B \\ & R_n \ar[ul]^{\tau_A}
\ar[ur]_{\tau_B} \\ }
$$

\medskip

{\em Culler Vogtmann Space of $F_n$} or {\em Outer Space of rank $n$}
is the set CV$_n$ of equivalence classes of marked metric graphs of
rank $n$.

\medskip

It is common to consider standard representative of a given class by
taking volume one graphs (here volume means total edge length.)

However, we usually do not normalise metric graphs, and when we will
do it we will use different normalisations depending on the
calculations we are making.

We note that since the equivalence allows homothety, given a point
$[X]$ in CV$_n$, we only have the metric on $X$ up to scaling
constants.  If one instead only considers the equivalence up to
isometry, then one obtains {\em unprojectivised} CV$_n$ and the metric
on the graph corresponding to a point there is determined by the
point.

\begin{rem} In the following, if there is no ambiguity, we will not
distinguish between a metric graph $X$ and its class $[X]$. If we need
to choose a particular representative of $[X]$ we will explicitly
declare that.
\end{rem}

\subsection{The Topology of CV$_n$} Outer space is endowed with
topology induced by edge-lengths of graphs.

Given any marked graph $A$, we can look at the universal cover $T_A$
which is an $\mathbb{R}$-tree on which $\pi_1(R_n)$ acts by
isometries, via the marking $\tau_A$.  Conversely, given any minimal
free action of $F_n$ by isometries on a simplicial $\mathbb{R}$-tree,
we can look at the quotient object, which will be a graph, $A$, and
produce a homotopy equivalence $\tau_A: R_n \to A$ via the
action. Equivalence of graphs in $CV_n$ corresponds to actions which
are equivalent up to equivariant homothety.

Thus, points in $CV_n$ can be thought of as equivalence classes of
minimal free isometric actions on simplicial $\mathbb{R}$-trees.
Given an element $w$ of $F_n$ and a point $A$ of the unprojectivised
$CV_n$, with universal cover $T_A$ whose metric we denote by $d_A$, we
may consider,
$$
L_A(w):=\inf_{p \in T_A} d_{A}(p, wp).
$$

It is well known that this infimum is always obtained and that, for a
free action, it is non-zero for the non-identity elements of the
group. In this context, $L_A(w)$ is called the translation length of
the element $w$ in the corresponding tree and clearly depends only on
the conjugacy class of $w$ in $F_n$.  If we look at graph $A$, then
$L_A(w)$ is the length of the geodesic representative of $w$ in $A$,
that is to say, the length of shortest closed loop representing free
homotopy class of ${\tau_A}_*(w)$ as an element of $\pi_1(A)$.  Thus
for any point, $A$, in $CV_n$ we can associate the sequence
$(L_A(w))_{w \in F_n}$ and it is clear that equivalent marked metric
graphs will produce two sequences, one of which is a multiple of the
other by a positive real number (the homothety constant). Moreover, it
is also the case that inequivalent points in $CV_n$ will produce
sequences which are not multiples of each other \cite{CuMo}. Thus, we
have an embedding of $CV_n$ into $\mathbb{R}^{F_n}/ \sim$, where
$\sim$ is the equivalence relation of homothety. The space $CV_n$ is
given the subspace topology induced by this embedding.

Finally it is clear we can realise any automorphism, $\phi$, of $F_n$
as a homotopy equivalence, also called $\phi$, of $R_n$. Thus the
automorphism group of $F_n$ acts on $CV_n$ by changing the marking.
That is, given a point $(A, \tau_A)$ of $CV_n$ the image of this point
under $\phi$ is $(A, \tau_A \phi)$.

$$
\xymatrix{ R_n \ar[r]^{\phi} \ar@/^2pc/[rr]|{\tau_A \phi} & R_n
\ar[r]^{\tau_A} & A.  \\ }
$$

Since two automorphisms which differ by an inner automorphism always
send equivalent points in $CV_n$ to equivalent points, we actually
have an action of $Out(F_n)$ on $CV_n$, and this space is called {\em
Outer Space} for this reason.

\subsection{Simplicial Subdivision of CV$_n$}

Given a rank-$n$, marked, metric graph $X$ whose edges are labelled
$e_1,\dots,e_k$, we can consider all marked metric graphs homeomorphic
to $X$ and with same marking. Such subset of CV$_n$ can be embedded in
$\mathbb R^k$ by
$$X\mapsto (L_X(e_1),\dots,L_X(e_k)).$$

If we consider standard normalisation with volume one, we obtain
standard open $(k-1)$-simplex of $\mathbb R^k$, i.e. the set
$\{(x_1,\dots,x_k)\in\mathbb R^k: x_i>0, \sum x_i=1\}$.

This gives us a natural subdivision of CV$_n$ into open simplices.

\begin{dfn} Let $\Delta$ be an open simplex of CV$_n$. The {\em
(marked) graph underlying} of $\Delta$ is the (marked) topological
type of graphs corresponding to points of $\Delta$.
\end{dfn}

Simplices of CV$_n$ will have some ideal faces and some true
faces. More precisely, in an abstract way, if $\Delta$ is a simplex
with underlying graph $X$, a face $\delta$ of $\Delta$ is obtained by
setting to zero the lengths of some of the edges of $X$. This
topologically corresponds to collapsing such edges. If the resulting
graph has still rank $n$, then $\delta$ exists as a simplex of CV$_n$,
and in this sense it is a true face. On the other hand, if the rank
decreases, then $\delta$ is not in CV$_n$ (and in fact belongs to the
boundary at infinity of CV$_n$) and in this case we say that $\delta$
is an ideal face of $\Delta$.

In what follows we always deal with true faces.

\begin{dfn} Let $\Delta$ be a simplex of CV$_n$ with underlying marked
graph $X$. A {\em face} of $\Delta$ is a simplex of CV$_n$ whose
underlying marked graph is obtained from $X$ by collapsing some
edges. The codimension of the face of $\Delta$ is the number of
collapsed edges.
\end{dfn}

It is readily checked by an Euler characteristic count that simplices
of maximal dimension of CV$_n$ correspond to trivalent graphs, and
that such graphs have $3n-3$ edges and $2n-2$ vertices. Therefore
their dimension is $3n-4$ as CV$_n$ is the projectivised outer
space. Looking at the topology of graphs we see that in general,

\begin{lemma}\label{l1dim} $k$-dimensional simplices of
(projectivised) CV$_n$ correspond to graphs with $k+1$ edges and
$k-n+2$ vertices.
\end{lemma}

Next we consider the $i$-skeleton of CV$_n$.

\begin{dfn} For $i\leq 3n-n$, the {\em $i$-skeleton CV$^i_n$} of
CV$_n$ is the set of simplices of CV$_n$ of dimension at most $i$.
\end{dfn}

An easy but important fact is that $i$-simplices correspond to smooth
points of the $i$-skeleton.

\begin{dfn} A point $x\in$CV$^i_n$ is {\em smooth} if it has a
neighbourhood in CV$^i_n$ homeomorphic to $\mathbb R^i$.
\end{dfn}

\begin{lemma}\label{l1sp} Open $i$-simplices of CV$_n$ are exactly the
connected components of the set of smooth points of CV$^i_n$. That is
to say
$$\{x\in CV_n^i:x \textrm{ is smooth } \}=\bigsqcup_{\Delta \textrm{
    open } i\textrm{-simplex}} \Delta
$$
\end{lemma} \proof It is enough to show that any $i-1$ simplex is the
face of at least three different $i$-simplices. Let $X$ be a point of
an $(i-1)$-simplex. Then $X$ is obtained by collapsing to zero an edge
$e$ of a point $\bar X$ of an $i$-simplex. Let $v_-$ and $v_+$ be the
endpoints of $e$. Clearly $v_-\neq v_+$ because otherwise the collapse
would decrease the rank. By definition, both $v_-$ and $v_+$ have
valence at least three, and they are identified in $X$ to the same
vertex $v$ which therefore has valence at least four.

For any subdivision of the set of germs of edges at $v$ in two subsets
of at least two germs, we can form a different $i$-simplex, having $X$
in one of its faces, by separating such subsets and inserting a new
edge between them. Clearly, different subdivisions give different
$i$-simplices, and we have at least three such subdivisions because
the valence of $v$ is at least four.\qed

\subsection{Roses and Multi-thetas}

Our result will be based on a detailed study of isometries of two
particular classes of marked graphs. Namely roses and multi-theta
graphs.

\begin{dfn} A {\em rose simplex} is a simplex $\Delta$ of CV$_n$ whose
underlying graph is a rose, i.e. a bouquet of $n$ copies of $S^1$.
Edges of such a graph are also called {\em petals}.  The {\em centre}
of $\Delta$ is the symmetric graph, that is to say one whose petals
all have the same length.
\end{dfn}

One should note that in the definition above, the centre is only
defined by specifying that the edges have the same length without
saying what that length is. We shall usually take a representative
whose petals all have length 1 but the reader should be aware that as
long as all the petals have the same length, the point in CV$_n$ will
be the same.

By Lemma~\ref{l1dim}, rose simplices are those simplices of lowest
dimension of CV$_n$.

\begin{dfn}
\label{centre} A {\em multi-theta simplex} is a simplex $\Delta$ of
CV$_n$ whose underlying graph has only two vertices and $n+1$ edges
joining them (such graph is called a multi-theta.) The centre of
$\Delta$ is the symmetric graph, that is to say, the one whose edges
have all same length.
\end{dfn}

\setlength{\unitlength}{.8ex}
\begin{figure}[htbp] \centering
  \begin{picture}(30,20) \qbezier(0,10)(15,30)(30,10)
\qbezier(0,10)(15,20)(30,10) \qbezier(0,10)(15,10)(30,10)
\qbezier(0,10)(15,0)(30,10) \qbezier(0,10)(15,-10)(30,10)
  \end{picture}
  \caption{A multi-theta graph in CV$_4$}
\end{figure}
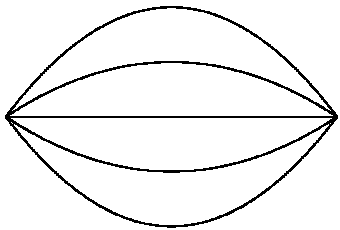

\begin{dfn} A rose-face of a simplex $\Delta$ of CV$_n$ is a rose
simplex which is a face of $\Delta$.
\end{dfn}

Formally speaking, simplices are open, so the rose-face of a simplex
is not subset of it. Nonetheless, it is readily checked that any
isometry of a simplex extends to its faces and rose-faces, though it
may permute them. As we are interested in studying isometries, by
abuse of notation, we will consider the rose-faces of a simplex as
subsets of it.

\begin{rem} Let $\Delta$ be a simplex of CV$_n$ with underlying graph
$X$. Then any rose-face of $\Delta$ is obtained by collapsing a
maximal tree $T$ of $X$, and different trees give rise to different
faces.  Therefore rose-faces of $\Delta$ are in correspondence with
maximal trees of $X$ (for instance, in case of multi-theta simplices,
rose-faces are in correspondence with edges).
\end{rem}

\subsection{Distances and stretching factors} We recall here the
definitions of --- both the symmetric and non-symmetric --- Lipschitz
distances on CV$_n$. These are defined via stretching factors of maps
between points of outer space. Stretching factors, outer space and
related topics are widely studied by many authors, and literature on
the matter is huge (see for
instance~\cite{MR1147956,MR2299451,MR2197815,MR1182503,CuVo86,FrMa,MR2299450,MR2216713,MR2439402,MR2322181,MR2373015}.)

\begin{dfn} For any two points $X$ and $Y$ in CV$_n$, normalised to
have volume one, we define the right stretching factor
as $$\Lambda_R(X,Y)=\sup_\gamma\frac{L_X(\gamma)}{L_Y(\gamma)}$$ where
the supremum is taken over all loops (or, equivalently over all
conjugacy classes in $F_n$.)  Similarly, the left stretching factor is
$$\Lambda_L(X,Y)=\Lambda_R(Y,X)=\sup_\gamma\frac{L_Y(\gamma)}{L_X(\gamma)}.$$
\end{dfn}

\begin{dfn} For any two points $X$ and $Y$ in CV$_n$, normalised to
have volume one, the right and left distances are defined by
$$d_R(X,Y)=\log(\Lambda_R(X,Y))\qquad d_L(X,Y)=\log(\Lambda_L(X,Y)).$$
\end{dfn}

\begin{dfn} For any two points $X,Y\in$ CV$_n$, not necessarily
normalised, the symmetric bi-Lipschitz metric, is defined by
$$d(X,Y)=d_R(X,Y)+d_L(X,Y)=\log
\sup_\gamma\frac{L_X(\gamma)}{L_Y(\gamma)}
\sup_\gamma\frac{L_Y(\gamma)}{L_X(\gamma)}.$$
\end{dfn}

We refer to~\cite{FrMa} for a detailed discussion on such metrics. We
recall some basic facts. Firstly, the suprema in definitions are
actually maxima.  Also, we recall that $Out(F_n)$ acts faithfully by
isometries on CV$_n$ (for $n\geq 3$, in rank two the kernel of the
action is $\mathbb Z_2$) endowed with any of above metrics. Finally we
note that the symmetric metric is scale invariant, while the
non-symmetric ones require normalisation.

The main tool for studying such distances is the so-called sausages
lemma, which allows us to quickly compute stretching factors, and
which we will use extensively throughout the paper (see~\cite{FrMa}
for the proof).

\begin{dfn}[Almost simple closed curves] Let $X$ be a point of
CV$_n$. A {\em simple closed curve} (s.c.c. for short) is an embedding
of $S^1$ to $X$. A {\em figure-eight} curve is an embedding to $X$ of
the bouquet $S^1\wedge S^1$ of two circles. A {\em barbell} curve is
roughly speaking an embedding to $X$ of the space: O---O. More
precisely, let $Q=\{(x,y)\in\mathbb R^2:\ \sup(|x|,|y|)=1\}$, then a
barbell curve is an immersion $c:Q\to X$ such that $c(x,y)=c(x',y')$
if and only if $x=x'$ and $|y|=|y'|=1$.

An {\em almost simple closed curve} (a.s.c.c. for short) is a curve
which is either an s.c.c., or a figure-eight or a barbell curve.
\end{dfn}

\begin{lemma}[Sausages Lemma]\label{lsausages} For any two marked
metric graphs $X$ and $Y$
$$\sup_\gamma\frac{L_Y(\gamma)}{L_X(\gamma)}$$ is realised by an
a.s.c.c. of $X$. Moreover, If both $X$ and $Y$ are roses, then the
supremum is realised by petals.
\end{lemma}

We notice that the Sausages Lemma not only allows to actually compute
distances, but is also important from a theoretical
view-point. Indeed, the fact that lengths of a.s.c.c. determine
distances, and therefore points of outer space, is a key-point in the
proof of Theorem~\ref{tmain} (see in particular
Theorems~\ref{tasymptotic} and~\ref{tsame}).

\medskip

Another simple but somehow surprising result that we will need in the
sequel is the following (whose proof can be found in~\cite{FrMa}).

\begin{lemma}\label{lgeod} Suppose $\sigma$ is a $d$-geodesic between
two points $X$ and $Y$ of CV$_n$. Let $Z$ be a point in $\sigma$. A
loop $\gamma_0$ is maximally (resp. minimally) stretched from $X$ to
$Y$ --- that is to say, it realises $\sup_\gamma
L_Y(\gamma)/L_X(\gamma)$ --- if and only if the same is true from $X$
to $Z$ and from $Z$ to $Y$.
\end{lemma}

\section{Schema of proof of Theorem~\ref{tmain}}\label{sschema} We
briefly describe here the strategy for proving our main result.  We
recall that we aim to show that any isometry $\Phi$ of CV$_n$ is
induced by some element of Out$(F_n)$.

\begin{enumerate}
\item For topological reasons, $\Phi$ maps simplices to
simplices. Moreover it maps rose simplices to rose simplices and
multi-theta simplices to multi-theta simplices.
\item Computation of isometry group of rose simplices (it will be
$\mathbb R^n\rtimes$ a finite group.)

\item For a point $X$ in a simplex $\Delta$ of CV$_n$, the asymptotic
behaviour of distances from $X$ to points in rose-faces of $\Delta$
determine lengths of simple closed curves of $X$. This being true not
only for points of $\Delta$ but also for points in any other simplex
having the same rose-faces as $\Delta$.
\item \label{four} For a point $X$ in a simplex $\Delta$ (or in other
simplices sharing rose-faces with $\Delta$) the lengths of simple
closed curves and the asymptotic behaviour of distances from $X$ to
points in rose-faces of $\Delta$, determine lengths of almost simple
closed curves of $X$ (whence asymptotic distances determine lengths of
a.s.c.c.)
\item Study of isometries of multi-theta simplices. We show that any
isometry of a multi-theta simplex fixes its centre. How:
  \begin{enumerate}
  \item Study of those pairs of points joined by a unique geodesic,
showing that for any point $X$ in the interior of $\Delta$, there is a
standard set of ``rigid'' geodesics emanating from $X$.
\item Show that for any point other than the centre, there is at least
one more ``rigid'' geodesic, while for the centre, the
standard set is all we have. This characterises the centre of $\Delta$
from a metric point of view.
\item Finally, for the centre of any rose-face of $\Delta$ there is a
unique ``rigid'' geodesic joining it to the centre.
  \item In particular, any isometry fixes the centre, and if it does
not permute rose faces, it fixes also such ``rigid''
geodesics.
  \end{enumerate}
\item Combining this with the knowledge of isometries of roses, we get
that if an isometry of a multi-theta simplex, $\Delta$, does not
permute its rose-faces, then it point-wise fixes them and hence
point-wise fixes $\Delta$ by~\ref{four} above (and we always can
reduce to the situation where $\Phi$ leaves some multi-theta $\Delta$
and all its rose faces invariant, by composing with an appropriate
element of $Out(F_n)$).

\item Show that any isometry that fixes a multi-theta simplex, also
fixes all rose simplices of CV$_n$ (not only its faces.)
\item Show that simplices that are possibly permuted by $\Phi$ share
their rose-faces, and that simplices that share rose-faces ``have the
same set of simple closed curves and the same set of almost simple
closed curves''. As asymptotic distances from rose-faces determine
points in such simplices, it follows a posteriori that they cannot be
permuted.
\end{enumerate}

\section{Topological constraints for homeomorphisms}

In this section we prove first step of our strategy, that isometries
of CV$_n$ respect its simplicial and incidence structure. That result
does not require any metric structure, just the fact that isometries
are homeomorphisms.

\begin{lemma}\label{L1} Any homeomorphism of CV$_n$ maps
$k$-dimensional simplices to $k$-dimensional simplices.
\end{lemma} \proof The proof goes by induction on the
codimension. Open top dimensional simplices coincide with smooth
points (Lemma~\ref{l1sp}.)

Clearly, to be a smooth point is invariant under
homeomorphisms. Again, by Lemma~\ref{l1sp}, open top-dimensional
simplices are exactly connected components of set of smooth points.
Therefore homeomorphisms map open top-dimensional simplices to open
top-dimensional simplices.

Suppose the claim true for dimensions greater than $i$.  By induction,
any homeomorphism $\Phi$ of CV$_n$ induces a homeomorphism of
$i$-skeleton CV$_n^i$. Open codimension-$(n-i)$ simplices are now
connected components of smooth part of CV$_n^i$, and therefore $\Phi$
maps $i$-simplices to $i$-simplices.\qed

\begin{lemma}\label{L2} Any homeomorphism of CV$_n$ maps
rose-simplices to rose-simplices, and multi-theta simplices to
multi-theta simplices.
\end{lemma} \proof This is just a dimensional argument. Clearly,
homeomorphisms preserve dimension. By Lemma~\ref{l1dim},
$n-1$-dimensional simplices are exactly rose-simplices, and the first
claim follows. If we look at $n$ dimensional simplices, we see that
multi-theta simplices are characterised by having exactly $n+1$
rose-faces. So the homeomorphic image of a multi-theta simplex still
is a multi-theta simplex.\qed

\section{Isometries of roses}\label{sisorose} In this section, we
compute the isometry groups of rose simplices. In rank two, it is
immediate to see that a rose simplex is isometric to $\mathbb R$, so
its isometries are known. For the general case we prove,

\begin{thm}\label{tisorose} The isometry group of a rose-simplex $R$
of CV$_{n+1}$ is $\mathbb R^n\rtimes \mathfrak F$, where $\mathfrak F$
is finite and stabilises the centre, and $\mathbb R^n$ acts
transitively.  Moreover, for $n\geq 2$, the group $\mathfrak F$ is $
S_{n+1}\times\mathbb Z_2$, where $S_{n+1}$ is the symmetric group on
$n+1$ letters and is induced by permutations of petals. For $n=1$
(i.e. in the rank-two case) $\mathfrak F=\mathbb Z_2= S_2$.
\end{thm} \proof

Any point of $R$ is determined by the lengths of its petals, that we
label $e_0,\dots,e_n$.  We identify the unprojectivised $R$ with
$\mathbb R^{n+1}$ as follows. To any $(x_0,\dots,x_n)\in\mathbb
R^{n+1}$ is associated the graph $X$ such that
$$L_X(e_i)=e^{x_i}$$

Note that origin of $\mathbb R^{n+1}$ corresponds to centre of
$R$. Moreover, scaling-equivalence on CV$_{n+1}$ descends to
relation $$x\sim y \qquad \textrm{ if and only if } \qquad
x-y=\lambda(1,\dots,1).$$

The pull back of the (pseudo) metric $d$ to $\mathbb R^{n+1}$ is then
$$d\big((x_0,\dots,x_n),(y_0,\dots,y_n)\big)=
\sup_i(x_i-y_i)+\sup_i(y_i-x_i).$$

This immediately implies that translations of $\mathbb R^{n+1}$ are
isometries, and that translations along vector $(1,\dots,1)$ are in
fact the only ones inducing the identity of the projectivised $R$. So
we have that
$$\mathbb
R^n=\mathbb R^{n+1}/<(1,\dots,1)>$$ acts freely and transitively on
$R$.

Thus, it remains to determine the stabiliser of the origin.

\medskip

Clearly, permutations of coordinates are isometries that fix
origin. Finally, we have the reflection
$$\sigma : (x_0,\dots,x_n)\mapsto (-x_0,\dots,-x_n).$$

In the rank-two case, that is to say when $n=1$, we are studying
isometries of $\mathbb R$ that fix origin. Therefore in rank-two, the
stabiliser of the origin consists of the reflection about the origin
and the identity: note that this reflection (the map $\sigma$, above,
in other words) is induced by the map which interchanges the two
petals of our rank 2 rose.

\medskip

For $n>1$, our claim is that the stabiliser of origin is
$$\mathfrak F=S_{n+1}\times <\sigma>.$$

For that, we need some work. First of all, note that the (pseudo)
metric $d$ on $\mathbb R^{n+1}$ is induced by the (pseudo) norm
$$||x||=d(0,x).$$

In order to make $||\cdot||$ a norm and $d$ a metric, for any point
$x\in\mathbb R^{n+1}$ we choose the $\sim$-representative of
$x+\mathbb R(1,\dots,1)$ that has $0$ as the first coordinate. We can
do that because $(x_0,\dots,x_n)\sim
(x_0,\dots,x_n)-x_0(1,\dots,1)$. This gives an isometry between $R$
and $\mathbb R^n$ with the following metric (still denoted by $d$)
$$d\big((x_1,\dots,x_n),(y_1,\dots,y_n)\big):=
d\big((0,x_1,\dots,x_n),(0,y_1\dots,y_n)\big).$$

We give now a more explicit description of that metric.

\begin{lemma}\label{ll1} For any set $I\subseteq\{1,\dots,n\}$ let
$R^I$ be the sector of $\mathbb R^n$ such that either $x_i\geq 0$ for
all $i\in I$ and $x_i\leq 0$ for all $i\notin I$, or vice versa.
Then, for $x\in R^I$
$$||x||=||x||_{\infty,I}+||x||_{\infty,I^c}$$
where $||x||_{\infty,I}=\sup_{i\in I}|x_i|$ and $I^c$ is the
complement of $I$ in $\{1,\dots,n\}$.
\end{lemma} \proof This is a straightforward calculation. Indeed, by
definition
$$||x||=\sup\{0,\sup_{i=1,\dots, n} x_i\}+\sup\{0,\sup_{i=1,\dots, n} -x_i\}$$
and, when $x\in R^I$, that equals
$||x||_{\infty,I}+||x||_{\infty,I^c}$.  \qed

\medskip

Our next step is an idea that we will return to throughout the paper,
and it is that the ``unique'' geodesics are rather rare and allow one
to determine the possible isometries.

\begin{rem}\label{r1} Note that $l^1$-norms naturally present
phenomena of non-uniqueness of geodesics. Namely, consider two
geodesic spaces $(X_1,d_1)$ and $(X_2,d_2)$, and their cartesian
product equipped with the sum metric
$d((x_1,x_2),(y_1,y_2))=d_1(x_1,y_1)+d_2(x_2,y_2)$. Then any geodesic
$\gamma:[0,1]\to X_1\times X_2$ is of the form
$\gamma=(\gamma_1,\gamma_2)$, and, up to reparametrisation,
$$t\mapsto\left\{\begin{array}{ll}(\gamma_1(t),\gamma_2(0))&
    t\in[0,1]\\ (\gamma_1(1),\gamma_2(t-1)& t\in [1,2]
\end{array}\right.  \quad
t\mapsto\left\{\begin{array}{ll}(\gamma_1(0),\gamma_2(t))& t\in[0,1]\\
(\gamma_1(t-1),\gamma_2(t)& t\in [1,2]
\end{array}\right.
$$
are two different geodesics whenever neither $\gamma_1$ nor $\gamma_2$
is the constant map. This situation is exactly the one arising in each
sector $R^I$ as above, where, by Lemma~\ref{ll1}, we have the sum of
two $l^\infty$-norms.

\end{rem}

\begin{prop}
\label{ungeod} Let $x=(x_1, \ldots, x_n) \in \mathbb{R}^n$, equipped
with the metric $d$ above.  Then there exists a unique geodesic
joining the origin to $x$ if and only if there exists a real number
$\lambda$ such that for all $i$, $x_i = \lambda$ or $x_i=0$. This
geodesic is given (up to reparametrisation) by the path $\gamma_x$
whose $i^{th}$ coordinate at time $t$ is $tx_i$.

Equivalently, a point $x=(x_0, x_1, \ldots, x_n) \in \mathbb{R}^{n+1}$
represents a point in $\mathbb{R}^n$ joined to the origin by a unique
geodesic if and only if there exist $\lambda, \mu$ such that each
$x_i$ is equal to either $\lambda$ or $\mu$.
\end{prop}

\proof The last statement follows trivially from the first, on taking
the representative with $x_0=0$, obtained by subtracting one of
$\lambda$ or $\mu$ from each coordinate.

Next, let $O$ denote the origin of $\mathbb R^n$. For any $x,y$ let
$\overline{xy}$ denote the path whose $i^{th}$ coordinate at time $t$
is $x_i+t(y_i-x_i)$, $t\in[0,1]$.  By Remark~\ref{r1}, if there is a
set of indices $I$ such that $||x||_{\infty,I}||x||_{\infty,I^c}\neq
0$ then $x$ is joined to $O$ by at least two different geodesics.
Thus, up to rearranging coordinates and possibly applying the isometry
$\sigma$ above, we can suppose $0\leq x_1\leq\dots\leq x_n$. Clearly,
$d(\gamma_x(s),\gamma_x(t))=x_n|t-s|$, so that
$\gamma_x=\overline{Ox}$ is a geodesic.  Suppose there is $i$ such
that $0<x_i<x_n$. Then, consider the point
$x_\varepsilon=(x_1/2,\dots,x_i/2+\varepsilon,\dots,x_n/2)$. For small
enough $\varepsilon$ the path $\gamma_\varepsilon$ resulting on the
union of $\overline{Ox_\varepsilon}$ and $\overline{x_\varepsilon x}$
is a geodesic from $x$ to $O$ as
$d(\gamma_\varepsilon(s),\gamma_\varepsilon(t))=x_n|t-s|/2$. Also,
$\gamma_\varepsilon$ is not a reparametrisation of $\gamma_x$ because
they differ in their middle points.

Conversely, suppose that there is $i$ so that $x_j=0$ for $j<i$ and
$x_j=x_n$ for $j\geq i$. Let $\gamma$ be a geodesic between $O$ and
$x$. If there is a time $t$ such that the $j^{th}$ coordinate of
$\gamma(t)$ is different from $0$ for some $j<i$, then a direct
calculation shows that $d(0,\gamma(t))+d(\gamma(t),x)$ is strictly
bigger than $x_n$ (while $d(O,x)=x_n$.) Thus, the $j^{th}$ coordinates
of $\gamma(t)$ all vanish for $j<i$. The very same argument shows that
for $j\geq t$ the $j^{th}$ coordinate of $\gamma(t)$ equals the
$n^{th}$ one, this showing that $\gamma_x$ is the unique geodesic from
$0$ to $x$.\qed

\medskip

Proposition~\ref{ungeod} is a translation of the fact that two roses
in the same simplex are joined by a unique geodesic if and only if
there are only two possible stretching factors for petals.

We note that Proposition~\ref{ungeod} gives us a collection of
geodesics which are permuted by any isometry fixing the origin. Using
this fact, we now proceed to calculate the stabiliser of the
origin. Since we already have that these geodesics must be permuted by
any isometry fixing the origin, we shall proceed by studying points on
these geodesics at fixed distance 1 from the origin. These are also
permuted and will give us the information we need about the
stabiliser.

For any $I\subseteq \{1,\dots,n\}$ we define points $p_I^+$ and
$p_I^-$ in $\mathbb R^n$ by
$$p_I^\pm=(x_1,\dots,x_n):\qquad x_i=\left\{
  \begin{array}{ll} \pm 1 & i\in I\\ 0 & i\in I^c
  \end{array}\right.
$$

such points are equivalents to points $P_I$ of $\mathbb R^{n+1}$
$$P_I=(x_0,\dots,x_n):\qquad x_i=\left\{
  \begin{array}{ll} 1 & i\in I\\ 0 & i\in I^c
  \end{array}\right.
$$

where $p^+_I$ is equivalent to $P_{0\cup I}$, and $p^-_I$ is
equivalent to $P_{I^c}$ (the complement here is made in
$\{0,\dots,n\}$.)

\begin{lemma}\label{lPI} For any distinct $I,J\subseteq \{0,\dots,n\}$
we have
$$d(P_I,P_J)=\left\{
  \begin{array}{ll} 1 & \textrm{ if }I\subseteq J \textrm{ or }
J\subseteq I\\ 2 & \textrm{ otherwise }
  \end{array}\right.
$$
Moreover, if $d(P_I,P_J)=2$, then the points of $\mathbb R^n$
corresponding to $P_I$ and $P_J$ are joined by a unique geodesic if
and only and $I=J^c$.
\end{lemma}

\proof The first part is a simple calculation. For the second part, we
use the fact that translations are isometries. Translate the point
$P_I$ to the origin and look at the image of $P_J$, which we call
denote $x=(x_0, \ldots, x_n)$. Then there will be a unique geodesic
between $P_I$ and $P_J$ if and only if there is a unique geodesic
between $x$ and the origin. However, is clear what each $x_i$ will
be. Namely,
$$
x_i = \left\{ \begin{array}{rl} 0 & \mbox{\rm if } i \in I \cap J \\ 0
& \mbox{\rm if } i \in I^c \cap J^c \\ 1 & \mbox{\rm if } i \in I^c
\cap J \\ -1& \mbox{\rm if } i \in I \cap J^c \\
\end{array} \right.
$$
As $d(P_I,P_J)=2$, we cannot have either $I\subseteq J \textrm{ or }
J\subseteq I$ and hence both $1$ and $-1$ must be taken by some of the
$x_i$. So by Proposition~\ref{ungeod}, $P_I$ and $P_J$ will be joined
by a unique geodesic if and only if no $x_i$ is equal to zero, which
is the same as saying $I \cap J = \emptyset = I^c \cap
J^c$. Equivalently, $I=J^c$.  \qed

\medskip

As stated, by Proposition~\ref{ungeod} any isometry that fixes the
origin must permute the $P_I$'s. For such an isometry $F$ and
$I\subset\{0,\dots,n\}$, we denote by $F(I)$ the set corresponding to
point $F(P_I)$.

From Lemma~\ref{lPI} we get
\begin{equation}\label{eIJ} \big(I\subseteq J \textrm{ or } J\subseteq
I\Big) \Leftrightarrow \Big(F(I)\subseteq F(J) \textrm{ or }
F(J)\subseteq F(I)\Big)
\end{equation} and
\begin{equation}\label{ec} F(I^c)=F(I)^c
\end{equation}

\begin{rem} The isometry $\sigma$ corresponds to $I\mapsto I^c$.
\end{rem}

\begin{lemma}\label{lcard} For any isometry $F$, the cardinality
$|F(I)|$ is either $|I|$ or $n+1 -|I|$.
\end{lemma} \proof

By $(\ref{eIJ})$ sets $I$ and $F(I)$ must have the same numbers of
subsets and supersets. For $I$ such number is $2^{|I|}+2^{n+1-|I|}-1$,
whence
$$2^{|I|}+2^{n+1-|I|}=2^{|F(I)|}+2^{n+1-|F(I)|}.$$

Set $x=\min\{|I|,n+1-|I|\}$ and $y=\min\{|F(I)|,n+1-|F(I)|\}$.  We
have
$$2^x(1+2^{k})=2^y(1+2^h)$$
for some non-negative numbers $k,h$. Whence $x=y$ and the claim
follows.\qed

\medskip

\begin{rem}
\label{invert} Up to possibly composing with $\sigma$ we may suppose,
as we do, that there is $i_0$ such that $|F(\{i_0\})|=1$.
\end{rem}

\begin{lemma}
\label{allone} If there is $i_0$ such that $|F(\{i_0\})|=1$, then for
all $i$ we have that $|F(\{i\})|=1$.
\end{lemma} \proof Note that by $(\ref{ec})$, $|F(\{i_0\}^c)|=
|F(\{i_0\})^c|= n$.  Now consider some $i\neq i_0$, whence
$\{i\}\subseteq \{i_0\}^c$. If $\{i\}=\{i_0\}^c$ then $n=1$ and the
lemma is proved. So we can suppose $\{i\}\neq\{i_0\}^c$, so
$F(\{i\})\neq F(\{i_0\}^c)=F(\{i_0\})^c$ (latter equality is by
$(\ref{ec})$.) Thus, by $(\ref{eIJ})$ and Lemma~\ref{lcard} we have
that $F(\{i\})$ is strictly contained in $F(\{i_0\})^c$. We therefore
have $|F(\{i\})|\leq n-1$, which implies $|F(\{i\})|=1$ because of
Lemma~\ref{lcard}. \qed

\medskip

\begin{rem} When $|F(\{i\})|=1$ for all $i$, we can define an element
$f$ of $ S_{n+1}$ by
$$F(\{i\})=\{f(i)\}.$$
\end{rem}

We show now that the permutation $F$ is actually induced by $f$.
\begin{lemma}
\label{perm} Suppose $|F(\{i\})|=1$ for all $i$. For all
$I\subseteq\{0,\dots,n\}$ we have
$$F(I)=\{f(i): i\in I\}.$$
\end{lemma} \proof For any $i\in I$ we have that $\{f(i)\}$ is either
contained in or contains $F(I)$, so we must have $f(i)\in F(I)$. The
same holds for $I^c$.\qed

\medskip

An immediate consequence of all these facts is the following fact.

\begin{prop}\label{pfix} Up to possibly composing with $\sigma$ and an
element of $S_{n+1}$, any isometry of $R$ that fixes origin also fixes
all points $p^\pm_I$.
\end{prop} \proof Let $F$ be an isometry of $R$ fixing the
origin. We shall also use $F$ to denote the induced permutation of
$\{0, \ldots, n \}$, so that $F(P_I)=P_{F(I)}$.

By remark~\ref{invert} and Lemma~\ref{allone}, we may suppose that for
all $i$ we have $|F(\{i\})|=1$. Hence by Lemma~\ref{perm}, $F$ is
induced by some permutation, $f$. We can think of this permutation as
an isometry of $R$ which permutes the petals of the rose. By composing
$F$ with the inverse of this isometry, we get that $F(I)=I$ for all
subsets $I$ of $\{0, \ldots, n \}$. Thus $F$ fixes all the points
$P_I$ and thus all the $p^\pm_I$. \qed

\medskip

Next lemma is a simple case of a general asymptotic argument (see
Section~\ref{sdr} and compare in particular with
Proposition~\ref{pasymptotic1}).

\begin{lemma}\label{lbus} For any $i=1,\dots, n$ and $t\in\mathbb R$,
let $x_i(t)$ be the point of rose $R$, identified with $\mathbb R^n$,
whose coordinates are zero except for $i^{th}$ which is $t$:
$$x_i(t)=(0,\dots,t,\dots,0), \qquad \textrm{t at the $i^{th}$ place}$$
and let $x_0(t)=(t,\dots,t)$.  Then, points of $R$ are determined by
distances from points $x_i(t)$'s.
\end{lemma} \proof Let $y=(y_1,\dots,y_n)$ be a point of $R$. Clearly,
for large enough $t$, we have
$$d(y,x_i(t))=t-y_i + \max\{0,\sup_{j\neq i} y_j\}.$$
Therefore, by knowing such distances, we know for each $i$
$$-y_i+\max\{0,\sup_{j\neq i} y_j\}.$$
Note that the $y_i$'s are all negative numbers if and only if such
quantities are all positive, and in that case they give exactly
$-y_i$.  On the other hand, if some non-positive quantity appears,
then indices $i$ for which $y_i$ is maximum are characterised by the
fact that $i^{th}$ quantity is not positive. Thus, varying $i>0$ we
know all differences $y_i-y_j$ for any $i,j$.

Finally, consider distances from $x_0(t)$ as $t\to-\infty$
$$d(y,x_0(t))=\max\{0,\sup_i(t-y_i)\}+\max\{0,\sup_i(y_i-t)\}=\max_i y_i-t.$$

This gives knowledge of $\max_i y_i$, and since we know those indices
for which $y_i$ is maximum, and all differences $y_i-y_j$, we get all
the $y_i$'s.  \qed \medskip

Note that such a result can be re-paraphrased by saying that Busemann
functions of ideal vertices determines points.

\medskip

We are now able to finish proof of Theorem~\ref{tisorose}.  Let $\phi$
be an isometry of $R$. Up to composing with a translation of $\mathbb
R^{n}$, we can suppose that $\phi$ fixes the origin.  By
Proposition~\ref{pfix} after possibly composing with elements of $
S_{n+1}\times\langle\sigma\rangle$, we can suppose that $\phi$ fixes
all the points $p^\pm_I$. Therefore, by Proposition~\ref{ungeod},
$\phi$ must fix all the points $x_i(t)$ and $x_0(t)$.
Lemma~\ref{lbus} now implies that $\phi$ is the identity.  \qed

\medskip

We conclude this part anticipating results of subsequent sections. We
have seen what the isometry group of a rose simplex is, and we have
seen in particular that there are isometries which are not induced by
elements of $Out(F_{n+1})$. This seems, a priori, to count as evidence
against our final result. However, no such isometry arises as the
restriction of a global isometry of CV$_{n+1}$. Indeed, we will show
that translations of $\mathbb R^n$ and reflection $\sigma$ are not
restrictions of global isometries. On the other hand, any isometry in
$S_{n+1}$ is induced by a permutation of generators and hence by an
element of $Out(F_{n+1})$ (see Sections~\ref{sdr} and~\ref{stheta}, in
particular Remark~\ref{rrn} and Lemma~\ref{lsigma}).

Note that this is not enough to show that Isom$($CV$_{n+1})$ is
$Out(F_{n+1})$. Indeed, it could be possible that an isometry permutes
simplices, and second, that restrictions of an isometry to different
simplices are restrictions of different elements of $Out(F_{n+1})$. We
will see that this is not the case (Section~\ref{stheta}).

\section{Asymptotic distances from roses and global
isometries}\label{sdr}

In this section we generalise calculations made in
Section~\ref{sisorose} about the asymptotic behaviour of distances.
The underlying philosophy is that Busemann functions of ideal vertices
are enough to distinguish points of outer space.

\medskip

What we have in mind is to prove the following fact, that if an
isometry fixes all rose-simplices of CV$_n$ then it must be the
identity. This, together with results of next section, opens the way
towards Theorem~\ref{tmain}.

\medskip

The first point in proving that result is that a priori, an isometry
that fixes all rose-faces, could permute other simplices.

For that, we have to understand simplices that are possibly not
invariant under the action of such an isometry. Lemma~\ref{lsr} below
will tell us that any two putatively permuted simplices must have the
same rose-faces.

Then, our aim will be to show that a point $X$ is determined by
asymptotic distances from points in the rose-faces of the simplex
containing $X$. More precisely, we show how such distances determine
the lengths, in $X$, of every almost simple closed curve.

We emphasise that the results we are proving here (out of necessity,
due to Lemma~\ref{lsr}) depend only on the set of rose-faces, and not
on the simplex containing $X$.

That is to say, suppose $X$ and $Y$ are points of simplices $\Delta_1$
and $\Delta_2$ who share their rose-faces.  If for any $p$ in any
rose-face we have $d(X,p)=d(Y,p)$, then we show that for any two
a.s.c.c. $\gamma_1$ and $\gamma_2$
$L_X(\gamma_1)/L_X(\gamma_2)=L_Y(\gamma_1)/L_Y(\gamma_2)$ (so lengths
of a.s.c.c.  are equal up to scaling.)  Of course, we need also to
show that whenever $\Delta_1$ and $\Delta_2$ share rose-faces, then a
loop $\gamma$ is a.s.c.c. in $\Delta_1$ if and only if the same
happens in $\Delta_2$.

Since the distance from $X$ to $Y$ is computed using only
a.s.c.c. (because of Lemma~\ref{lsausages}) we deduce that this
implies $X=Y$.

\medskip

We start by studying simplices possibly permuted by isometries that
fix roses.

\begin{lemma}\label{lsr} Let $\Phi$ be an isometry of CV$_n$ that
fixes all rose-simplices.  If $\Delta$ is any simplex of CV$_n$, then
$\Delta$ and $\Phi(\Delta)$ have the same rose-faces.
\end{lemma} \proof Let $R$ be a rose-face of $\Delta$, then
$d(R,\Delta)=0$, so
$$0=d(R,\Delta)=d(\Phi(R),\Phi(\Delta))=d(R,\Phi(\Delta)).$$

Thus, $R$ is a rose-face also of $\Phi(\Delta)$. Using $\Phi^{-1}$ we
get the converse.\qed

\medskip

Now, we show that two simplices that share rose-faces have the same
a.s.c.c.

\begin{thm}\label{tsame} Let $\Delta_1$ and $\Delta_2$ be two
simplices of CV$_n$ that share their rose-faces. Then they have the
same set of almost simple closed curves. More precisely, if $\gamma$
is a conjugacy-class in $F_n$, then its geodesic representative in
$\Delta_1$ is simple if and only if it is simple in $\Delta_2$, and it
is a figure-eight or bar-bell curve in $\Delta_1$ if and only if the
same is true in $\Delta_2$ (possibly bar-bells become figure-eight
curves and vice versa.)
\end{thm}

\proof Let $G_1$ and $G_2$ be marked graphs corresponding to simplices
$\Delta_1$ and $\Delta_2$.  Any rose-face of $\Delta_i$ is obtained by
collapsing a maximal tree in $G_i$.

We first prove that a loop is simple in $G_1$ if and only if it is
simple in $G_2$.  Let $\gamma$ be a simple loop in $G_1$, and let
$e_1$ be an edge of $\gamma$. As $e_1$ is part of a simple loop, it
does not disconnect $G_1$. Extend $\gamma\setminus e_1$ to a maximal
tree $T_1$ in $G_1$. Let $R$ be the rose obtained by collapsing
$T_1$. The class of $\gamma$ in $R$ is represented by a petal $p$ (the
image of $e_1$.)  As $\Delta_1$ and $\Delta_2$ share rose-faces, $R$
is obtained by collapsing a maximal tree $T_2$ in $G_2$. So the class
of $\gamma$ in $G_2$ is represented by an edge $e_2$ corresponding to
the petal $p$ plus a path in $T_2$. As $T_2$ is a tree, such path is
unique and its union with $e_2$ is simple. Thus, $\gamma$ is
represented by a simple loop also in $G_2$.

\medskip Now, we deal with figure-eight and barbell curves. Let
$\gamma$ be such a curve in $G_1$. Let $\alpha$ and $\beta$ be the two
simple loops of $\gamma$.

\begin{lemma}\label{lmark} Let $R$ be any rose-face of $G_1$, then
$\alpha\cup\beta$ is represented in $R$ by a union of petals, each
petal appearing at most once. In particular the representatives of
$\alpha$ and $\beta$ in $R$ have no common petal.
\end{lemma} \proof Let $T$ be the maximal tree of $G_1$ collapsed in
order to obtain $R$.  Since $T$ is a tree, and it is maximal, it
cannot contain the whole $\alpha$, nor the whole $\beta$. As $\alpha$
and $\beta$ have no common edge, their images in $R$ share no
petal. Moreover, since $\alpha$ and $\beta$ are simple, no petal can
occur twice.  \qed

\medskip Note that Lemma~\ref{lmark} would fail if $\alpha\cup\beta$
were a theta curve. We can now conclude proof of Theorem~\ref{tsame}.
Let $e_\alpha$ be and edge of $\alpha$ and $e_\beta$ be an edge of
$\beta$. As $e_\alpha$ and $e_\beta$ are part of simple loops with no
common edges, we have that $e_\alpha\cup e_\beta$ does not disconnect
$G_1$. Extend $\gamma\setminus (e_\alpha\cup e_\beta)$ to a maximal
tree $T_1$, and let $R$ be the rose obtained by collapsing $T_1$. Let
$T_2$ be the tree of $G_2$ whose collapsing gives $R$. The loop
$\alpha$ is represented in $G_2$ by an edge corresponding to
$e_\alpha$, which we still denote by $e_\alpha$, and a path
$\sigma_\alpha$ in $T_2$ joining the end-points of $e_\alpha$. The
same (with the same notation) for $\beta$. The paths $\sigma_\alpha$
and $\sigma_\beta$ have connected intersection because $T_2$ is a
tree. It follows that the representative of $\gamma$ in $G_2$ is
either a figure-eight or a barbell, or a theta-curve. We show now that
the case of theta-curve cannot arise.

Indeed, suppose representative of $\gamma$ in $G_2$ is a
theta-curve. This is equivalent to saying that
$\sigma_\alpha\cap\sigma_\beta$ contains at least one edge
$e_0$. Clearly, $e_0$ does not disconnect $G_2$. We can therefore find
a maximal tree $T_0$ not containing $e_0$. Collapsing $T_0$ we get a
rose $R$ with a petal $p_0$ corresponding to $e_0$. In $R$, loops
representing $\alpha$ and $\beta$ share petals $p_0$. By
Lemma~\ref{lmark} $R$ cannot be obtained from $G_1$, in contradiction
with hypothesis that $\Delta_1$ and $\Delta_2$ have same
rose-faces.\qed

\medskip

Our next goal is to show that asymptotic distances from rose-faces
determine points.  First, we show how to determine lengths of simple
closed curves using distances from rose-faces. After, we will deal
with a.s.c.c.

\begin{prop}[Distances from roses determine simple
loops]\label{pasymptotic1} Let $G_1$ and $G_2$ be the underlying
graphs of two simplices $\Delta_1$ and $\Delta_2$ having the same
rose-faces. Let $X_1, X_2 \in \Delta_1 \cup \Delta_2$ such that
$d(X_1, Y)=d(X_2, Y)$ for any point $Y$ of any rose face of $\Delta_1$
(or $\Delta_2$).

Now fix a conjugacy class $\gamma_0$ which is a simple loop in $G_1$
(and hence $G_2$) and suppose that $X_1, X_2$ are the representatives
for which $L_{X_1}(\gamma_0)=L_{X_2}(\gamma_0)=1$.

Then, for {\em any} conjugacy class $\gamma$ in $F_n$ which is
represented by a simple loop in $G_1$,
$$
L_{X_1}(\gamma)=L_{X_2}(\gamma).
$$
\end{prop} Recall that $\gamma$ is simple in $G_1$ if and only if it
is simple in $G_2$ because of Theorem~\ref{tsame}.
Proposition~\ref{pasymptotic1} will follow from next lemma.

\begin{lemma}\label{lasymptotic1} Let $R$ be a rose simplex in $CV_n$
and let $e, e_0$ be petals in the underlying graph of $R$. Set $Y_t$
to be the ray in $R$, consisting of roses in $R$ all of whose edges
except $e, e_0$ have length $1$, and such that at time $t$
$$L_{Y_t}(e)=t \qquad L_{Y_t}(e_0)=\frac{1}{t}.$$

Now consider an $X\in \Delta$, where $\Delta$ is a simplex of $CV_n$
whose underlying graph is $G$ and such that $R$ is a rose-face of
$\Delta$.

\medskip

Let $\gamma$ be the simple closed curve in $X$ corresponding to
$e_0$. Also let $\gamma_e$ be an a.s.c.c. in $X$ which minimises
$\frac{L_X(\gamma_e)}{n_e(\gamma_e)}$, where $n_e(\gamma_e)$ is the
number of times $\gamma_e$ crosses $e$ (when projected to $R$.)
 Then,
  \begin{equation}
    \label{eq:las}
L_X(\gamma)=C(e,X)\lim_{t\to\infty}\frac{e^{d(X,Y_t)}}{t^2}
  \end{equation}

 where $C(e,X)=\frac{L_X(\gamma_e)}{n_e(\gamma_e)}$ as above.
\end{lemma}
{\bf NOTE}: The ray $Y_t$ depends only on the edges $e, e_0$, and the
rose simplex $R$.

\proof

Let $T$ be the maximal tree in $G$ corresponding to the projection of
$\Delta$ to $R$. We can find lifts of the edges $e, e_0$ in $G$. We
continue to call these edges $e$ and $e_0$.

Now consider the ray $Y_t$. We let $t\to\infty$ and study the
asymptotic behaviour of $d(X,Y_t)$.

We claim that for sufficiently large $t$, the loop $\gamma$ is
maximally shrunk from $X$ to $Y_t$. Indeed, if $\sigma= e_1\dots e_k$
is a loop, then

\begin{equation}
  \label{eminmax} \frac{L_{Y_t}(\sigma)}{L_X(\sigma)}=
\frac{\displaystyle{\sum_{i: e_i=e_0}\frac{1}{t}+ \sum_{i: e_i=e}t+
\sum_{i: e_i\notin (T\cup e\cup e_0)}1} }{\sum_i L_X(e_i)}.
\end{equation}

Since $T$ is a tree, it cannot contain loops. Thus, if in $\sigma$
there is some $e_0\neq e_i\notin T$ the above stretching factor is
bounded below uniformly on $t$. On the other side, if $\sigma=\gamma$,
the stretching factor goes to zero as $t\to\infty$. Finally, if in
$\sigma$ there is no edge $e_i\notin T\cup e_0$, then $\sigma$ is a
multiple of $\gamma$ because $T$ is a tree, and so $\gamma$ is the
only way to obtain a simple loop from $e_0$ by adding edges of $T$.

Now, we look for maximally stretched loops.  As above, we compute
$L_{Y_t}(\sigma)/L_x(\sigma)$ for a generic a.s.c.c. $\sigma$ using
$(\ref{eminmax})$. If $\sigma$ does not contain $e$, then there is an
upper bound to the stretching factor and, as $t\to\infty$, it is
readily checked that if $(\ref{eminmax})$ is maximised, then for big
enough $t$ the ratio of $L_X(\sigma)$ over the number of occurrences
of $e$ in $\sigma$ is minimised, hence
$$
\frac{L_X(\sigma)}{n_e(\sigma)}=
\frac{L_X(\gamma_e)}{n_e(\gamma_e)}:=C(e,X).
$$

It follows, that for sufficiently large $t$, if $\sigma$ is maximally
stretched we have

\begin{equation}
  \label{edl} d(X,Y_t)=\log
\frac{L_{Y_t}(\sigma)}{L_X(\sigma)}\frac{L_X(\gamma)}{L_{Y_t}(\gamma)}
=\log\frac{L_X(\gamma)}{L_X(\sigma)}
t(n_e(\sigma)t+b+n_{e_0}(\sigma)\frac{1}{t})
\end{equation}

where $n_e(\sigma)$, $n_{e_0}(\sigma)$ are either $1$ or $2$ as
$\sigma$ is a.s.c.c., $b$ is the number of edges of $\sigma$ not
belonging to $T\cup e_0\cup e$. Whence,

$$\frac{L_X(\gamma)}{L_X(\sigma)}=\lim_{t\to\infty}\frac{e^{{d(X,Y_t)}}}{n_e(\sigma)t^2}
$$
so
$$
\begin{array}{rcc} L_X(\gamma)& = &
\frac{L_X(\sigma)}{n_e(\sigma)}\lim_{t\to\infty}\frac{e^{{d(X,Y_t)}}}{t^2}
\\ & = &
\frac{L_X(\gamma_e)}{n_e(\gamma_e)}\,\,\lim_{t\to\infty}\frac{e^{{d(X,Y_t)}}}{t^2}
\\ & = & C(e,X)\lim_{t\to\infty}\frac{e^{d(X,Y_t)}}{t^2}.
\end{array}
$$
and the lemma is proved.\qed

\proof[Proof of Proposition~\ref{pasymptotic1}]

Let $X$ be a point of either $\Delta_1$ or $\Delta_2$.  Let $\gamma$
and $\eta$ be two simple closed curves in $G_1$. Choose an edge $e_0$
in $\gamma$ (but not in $\eta$) and an edge $f_0$ in $\eta$ (but not
in $\gamma$). We can then find a maximal tree $T$ in $G$ which extends
$\gamma \cup \eta - (e_0 \cup f_0)$. Let $R$ be the corresponding rose
face of $\Delta$. Note that in any rose within this simplex $\gamma$
and $\eta$ each project to a single petal, which we will call $e_0$
and $f_0$ (these petals are also the projections of those edges).

Now assume that $n \geq 3$ so that we can find yet another petal, $e$,
distinct from $e_0, f_0$.

By Lemma~\ref{lasymptotic1}, there is a ray $Y_t$ such that,

$$
L_X(\gamma)=C(e,X)\lim_{t\to\infty}\frac{e^{d(X,Y_t)}}{t^2}.
$$

Similarly, there is a ray $Z_t$ such that,
$$
L_X(\eta)=C(e,X)\lim_{t\to\infty}\frac{e^{d(X,Z_t)}}{t^2}.
$$
Hence,
\begin{equation}
\label{limit}
\frac{L_X(\gamma)}{L_X(\eta)}=\lim_{t\to\infty}\frac{e^{d(X,Y_t)}}{e^{d(X,Z_t)}}.
\end{equation} Moreover, by Lemma~\ref{lasymptotic1}, this last
equation must hold for {\em any} $X$ which has $R$ as a rose face
(where we simply interpret $\gamma, \eta$ as conjugacy classes of
$F_n$) and thus certainly for any $X \in \Delta_1 \cup
\Delta_2$. Thus, for the $X_1, X_2$ in the statement of the
Proposition,

$$
\frac{L_{X_1}(\gamma)}{L_{X_1}(\eta)}=
\frac{L_{X_2}(\gamma)}{L_{X_2}(\eta)}
$$
for any two loops $\gamma, \eta$ which are simple in $G_1$ (and hence
$G_2$). Putting $\eta=\gamma_0$ proves Proposition~\ref{pasymptotic1}
when $n \geq 3$.

\medskip

Now consider the case $n=2$. Note that here, distinct simplices have
different collections of rose faces (so we need not worry about
$\Delta_2$). When the underlying graph of $X$ is a rose, there are
exactly two simple loops in $X$, $\gamma$ and $\eta$ and
Lemma~\ref{lasymptotic1} produces exactly two different rays $Y_t$
with limits, as in \ref{limit}, $\frac{L_X(\gamma)}{L_X(\eta)}$ and
$\frac{L_X(\eta)}{L_X(\gamma)}$ and the order of these is independent
of $X$. So Proposition~\ref{pasymptotic1} is true in this case.

Similarly, if the underlying graph of $X$ is a barbell, then there is
exactly one rose face and exactly two simple loops, $\gamma,
\eta$. Again, the limits from \ref{limit} will give
$\frac{L_X(\gamma)}{L_X(\eta)}$ and $\frac{L_X(\eta)}{L_X(\gamma)}$
and the lemma is again true in this case.

\medskip

Finally, if the underlying graph of $X$ is a theta curve, then $X$ has
exactly 3 edges, $x, y, z$, 3 rose faces and 3 simple closed curves,
$x \overline{y}, x\overline{z}, y\overline{z}$. There are then 6
possible rays as in Lemma~\ref{lasymptotic1}. However, each limit,
$$
\lim_{t\to\infty}\frac{e^{d(X,Y_t)}}{t^2},
$$
is equal to one of $\frac{L_X(x \overline{y})}{C(z, X)}, \frac{L_X(x
\overline{z})}{C(y, X)}, \frac{L_X(y \overline{z})}{C(x, X)}$. Also
note that $C(x, X)$ is simply the length of the shortest simple loop
in $X$ which crosses $x$, since an a.s.c.c. in $X$ is actually a
simple loop. Hence, $C(x, X)$ is equal to either $L_X( x
\overline{y})$ or $L_X( x \overline{z})$.

Thus, if $x \overline{y}$ is the shortest simple loop in $X$, then
$\frac{L_X(x \overline{y})}{C(z, X)}$ will be the smallest of the
three limits, $C(x, X)=C(y, X) = L_X( x \overline{y})$ and
conversely. From these observations, the Proposition follows
easily. Take $X_1, X_2$ with the same distances to rose faces. Then
the limits above, for $X_1, X_2$ respectively, produce the same
ordered results (however, the $C$ terms need to be evaluated in
different $X_i$'s).

Nevertheless, if without loss of generality, $x \overline{y}$ is the
shortest loop in $X_1$, then the limit $L_{X_1}( x \overline{y})/ C(z,
X_1)$ will be least, and thus so will $L_{X_2}( x \overline{y})/ C(z,
X_2)$ and hence $x \overline{y}$ must also be the shortest loop in
$X_2$. The Proposition now readily follows.  \qed

\begin{rem} We note that the constant $C$ depends on $X$ and on $e$,
but not on $\gamma$ or $e_0$.  Such a dependence is thus cancelled
when we consider the ratio $L_X(\gamma)/L_X(\eta)$, which therefore
actually depends only on asymptotic distances from $X$ to rose-faces.
\end{rem}

We now have sufficient tools for proving that asymptotic distances
from $X$ to rose-faces determine the lengths of all a.s.s.c., whence
determine $X$.

\begin{thm}\label{tasymptotic} Let $\Delta_1$ and $\Delta_2$ be
simplices of CV$_n$ with the same set of rose-faces. Let $G_1$ and
$G_2$ be the underlying graphs of $\Delta_1$ and $\Delta_2$
respectively. Let $\gamma_0$ be a simple loop in $G_1$ (whence its
representative in $G_2$ is a simple loop as well.) For any class $[X]$
of metric graphs in $\Delta_1\cup\Delta_2$ consider the representative
$X$ so that $L_X(\gamma_0)=1$. Now consider two such representatives,
$X_1, X_2 \in \Delta_1 \cup \Delta_2$ such that $d(X_1, Y ) = d( X_2,
Y) $ for any $Y$ in any rose face of $\Delta_1$. Then, for any
a.s.c.c. $\gamma$ in $G_1$ (and hence $G_2$),
$$
L_{X_1}(\gamma)=L_{X_2}(\gamma).
$$
\end{thm} \proof The proof is in the same spirit as
Proposition~\ref{pasymptotic1}, but the situation now it is a little
more complicated.

Let $X$ be a point in either $\Delta_1$ or $\Delta_2$. It will be
sufficient to show that we can calculate the length of any a.s.c.c. in
$X$ by only using distances to rose faces.

By Proposition~\ref{pasymptotic1}, we know that lengths in $X$ of
simple loops are determined via asymptotic distances to particular
sequences of points, not depending on $X$. Thus, we can suppose that
we already know the lengths of all simple loops in $X$, because we
have normalised so that $L_X(\gamma_0)=1$. Thus what remains is to
deal with figure-eight and barbell curves. Clearly, the length of a
figure-eight is determined via Proposition~\ref{pasymptotic1}. On the
other hand, Theorem~\ref{tsame} tells that a figure-eight curve in
$\Delta_1$ may become a barbell in $\Delta_2$. For this reason we
treat figure-eight and barbell curves at the same time, considering a
figure-eight as a barbell whose central segment is reduced to a point.

In order to do this, we proceed as in Proposition~\ref{pasymptotic1};
for any given barbell curve, we build an appropriate sequence of
points $Y_t$ in some rose-face, such that the asymptotic distances
from $Y_t$ determine the length of the barbell.

\begin{rem} At this point, the reader should be aware of the subtle
difference in the argument from that in
Proposition~\ref{pasymptotic1}. Indeed, the points $Y_t$ we
constructed in Lemma~\ref{lasymptotic1} do not depend on $X$, but just
on $e_0$ (hence on $\gamma$) and $e$. Here, the ray $Y_t$ we shall
define will actually depend on $X$, or at least seem to, and thus
present a logical obstacle to our argument.

  More precisely, the ray $Y_t$ here will depend on lengths of simple
loops in $X$. Intuitively speaking, the ray $Y_t$ escapes to infinity
in a rose face and the ``slope'' of the this ray is determined by the
lengths of simple loops in $X$. However, this is sound because of
Proposition~\ref{pasymptotic1}.  So for any barbell curve, the ray we
chose for computing its length is the same for both $X_1$ and $X_2$,
thus barbells have same lengths in $X_1$ and $X_2$, and
Theorem~\ref{tasymptotic} will be proved.
\end{rem}

The rank-two case is easy and left to the reader (just use the
following argument without the need to introduce the edge $e$ and the
loop $\gamma_e$.) Suppose $n\geq 3$.

Let $\gamma$ be a barbell curve, possibly degenerate to a figure-eight
curve, say in $G_1$. Let $\gamma_1$ and $\gamma_2$ be the two simple
loops of $\gamma$, and let $e_1\in\gamma_1$ and $e_2\in\gamma_2$ be
two edges. Clearly, $G_1\setminus (e_1\cup e_2)$ is connected. Extend
$\gamma\setminus (e_1\cup e_2)$ to a maximal tree $T$, and consider
the rose $R_T$ obtained by collapsing $T$. Since $n>2$ there is an
edge $e\notin (T\cup \gamma)$. Also, there is a simple loop not
containing $e$ (for instance, $\gamma_1$.)

In $R_T$ we still denote by $e,e_1,e_2$ the petals corresponding to
$e,e_1,e_2$ respectively.

Now, look at simplex $\Delta_2$. Since $R_T$ is a rose-face also of
$G_2$, it is obtained by collapsing a maximal tree $T'$ in
$G_2$. Therefore, petals $e,e_1,e_2$ correspond to edges of
$G_2\setminus T'$, and the representative of $\gamma$ in $G_2$ is
disjoint from $e$.

Now let $Y_t$ be the point of $R_T$ whose petals have length $1$
except $e,e_1,e_2$ for which we set
$$L_{Y_t}(e)=t\qquad L_{Y_t}(e_1)=\frac{L_X(\gamma_1)}{t}\qquad L_{Y_t}(e_2)=\frac{L_X(\gamma_2)}{t}.$$

We now let $t\to\infty$. If $\sigma= l_1\dots l_k$ is a loop, then
(replace $T$ with $T'$ if $X\in\Delta_2$)
\begin{equation}
  \label{eminmax2} \frac{L_{Y_t}(\sigma)}{L_X(\sigma)}=
\frac{\displaystyle{\sum_{i: l_i=e_1}\frac{L_X(\gamma_1)}{t}+ \sum_{i:
l_i=e_2}\frac{L_X(\gamma_2)}{t}+ \sum_{i: l_i=e}t+ \sum_{i: l_i\notin
(T\cup e_1\cup e_2)}1} }{\sum_i L_X(e_i)}.
\end{equation}

Note that by Lemma~\ref{lsausages}, the loop $\sigma$ minimising the
equation above is realised by an a.s.c.c. in $Y_t$ (note this
statement is independent of $t$) and inspection of the equation
\ref{eminmax2} shows that the only possible candidates are $e_1, e_2$
and $e_1 e_2$. It is then easy to see that $e_1 e_2$ is a loop which
realises the minimum, and this is exactly the realisation of $\gamma$
in $Y_t$. (Note that when the barbell is actually a figure-eight, all
three loops give the same answer, but our statement remains true.)

As in Lemma~\ref{pasymptotic1}, one also checks that, for large enough
$t$, any maximally stretched loop $\sigma$ from $X$ to $Y_t$ ({\em i.e.} one
that maximises \ref{eminmax2}) must minimise the ratio of
$L_X(\sigma)$ over the number of occurrences of $e$ in $\sigma$, among
all loops. Such a ratio is exactly the constant $C(e,X)$ introduced in
Lemma~\ref{pasymptotic1}.

Distances may then be computed, and we obtain an expression of the
form,

$$
d(X,Y_t)=\log
\frac{L_{Y_t}(\sigma)}{L_X(\sigma)}\frac{L_X(\gamma)}{L_{Y_t}(\gamma)}
=\log\frac{L_X(\gamma)}{L_X(\sigma)}
\frac{t(at+b+c\frac{1}{t})}{L_X(\gamma_1)+L_X(\gamma_2)}.
$$

where $a$ is the number of occurrences of $e$ in $\sigma$. Thus

$$\frac{L_X(\gamma)}{L_X(\gamma_1)+L_X(\gamma_2)}=
\frac{L_X(\sigma)}{a}\lim_{t\to\infty} \frac{e^{d(X,Y_t)}}{t^2}
=C(e,X)\lim_{t\to\infty} \frac{e^{d(X,Y_t)}}{t^2}
$$

Since $L_X(\gamma_1)$ and $L_X(\gamma_2)$, are known, we just need to
determine $C(e,X)$, which is given by Lemma~\ref{lasymptotic1} in
terms of asymptotic distances. Namely, if $(Z_t)$ is the sequence of
points given by Lemma~\ref{lasymptotic1} for computing the length of
$\gamma_1$ we get $L_X(\gamma_1)=C(e,X)\lim_t e^{d(X,Z_t)}/t^2$.

If one likes exact formulae, one would have to introduce sequences
$(Z^1_t)$ and $(Q^1_t)$, given by Lemma~\ref{lasymptotic1} for the
ratio $L_X(\gamma_1)/L_X(\gamma_0)$; then look at sequences $(Z^2_t)$
and $(Q^2_t)$, given by Proposition~\ref{pasymptotic1} for the ratio
$L_X(\gamma_2)/L_X(\gamma_0)$, and get (remembering the normalisation
$L_X(\gamma_0)=1$, and noting that the edge $e$ may occur in
$\gamma_0$ so that all the sequences below may be different)
$$
L_X(\gamma)=\lim_{t\to\infty}
\frac{e^{d(X,Y_t)}}{e^{d(X,Z_t)}}\frac{e^{d(X,Z^1_t)}}{e^{d(X,Q^1_t)}}(\frac{e^{d(X,Z^1_t)}}{e^{d(X,Q^1_t)}}+\frac{e^{d(X,Z^2_t)}}{e^{d(X,Q^2_t)}})
$$
\qed

\medskip

Finally, we are able to deal with global isometries of CV$_n$, proving
that isometries are determined by their restrictions to
rose-simplices.

\begin{thm}\label{tid} The only isometry of CV$_n$ that fixes all
rose-simplices is the identity.
\end{thm} \proof Let $\Phi$ be such an isometry. Let $X$ be a point of
a simplex $\Delta_1$ of CV$_n$, and let $\Delta_2=\Phi(\Delta_1)$. By
Lemma~\ref{lsr}, $\Delta_1$ and $\Delta_2$ share their rose-faces. By
Theorem~\ref{tsame} simple loops in $\Delta_1$ are also simple in
$\Delta_2$. In particular we can choose a simple loop $\gamma_0$ and
consider representatives of metric graphs of $\Delta_1$ and $\Delta_2$
by imposing that the length of $\gamma_0$ is $1$.

By Theorem~\ref{tsame}, $\Delta_1$ and $\Delta_2$ have the same almost
simple closed curves.  Since $\Phi$ fixes points in rose simplices,
for any $Y$ in a rose-face of $\Delta_1$, we have
$d(X,Y)=d(\Phi(X),\Phi(Y))=d(\Phi(X),Y)$.  Then,
Theorem~\ref{tasymptotic} says that the lengths of almost simple
closed curves are the same in $X$ and $\Phi(X)$. Therefore, the
Sausages Lemma~\ref{lsausages} implies $X=\Phi(X)$ (whence
$\Delta_1=\Delta_2$).\qed

\section{Isometries of multi-theta simplices and their
extensions}\label{stheta}

Recall that our main result is that isometries of Outer Space are all
induced by automorphisms of the free group. By Theorem~\ref{tid}, it
is enough to show that up to composing with automorphisms, we can
reduce to the case of isometries that point-wise fix every
rose-simplex, and we do that by studying isometries of multi-theta
simplices.

\medskip

Our first main result of this section is that isometries of
multi-theta simplices are induced by permutations of edges. Thus we
have no translations or inversions as in rose-simplices. In particular
this also shows that translations and inversions of rose-simplices
cannot arise as restrictions of global isometries of CV$_n$.

\medskip Then, we will prove that situation is in fact even more
rigid. Indeed, we show that if two isometries coincide on a
multi-theta simplex, then they coincide on all rose-simplices of
CV$_n$ (not only on faces of that simplex). This will basically
conclude Theorem~\ref{tmain}.

\medskip

We start by proving following theorem.

\begin{thm}\label{tmultitheta} Let $\Delta$ be a multi-theta simplex,
and let $\Phi$ be an isometry of $\Delta$.  Then $\Phi$ fixes the
centre of $\Delta$ (recall definition~\ref{centre}). Moreover, if
$\Phi$ leaves invariant all the rose-faces of $\Delta$, then it
actually fixes them point-wise, and in that case $\Phi$ is the
identity map on $\Delta$.
\end{thm}

Before proving Theorem~\ref{tmultitheta}, we need to establish some
preliminary technical lemmas. We follow the strategy sketched in
schema of Section~\ref{sschema}, focusing on the study of those pairs
of points that are joined by a unique geodesic.  We recall that Outer
Space is not a geodesic space; nevertheless, in any simplex, segments
(for the linear structure of the simplex) are geodesic. More
precisely, if we are in $CV_n$, then the points within a multi-theta
simplex $\Delta$ are specified by $n+1$ positive reals (giving an open
$n$-simplex, since one further needs to projectivise), corresponding
to the lengths of the $n+1$ edges. Then, given $x=(x_1, \ldots,
x_{n+1})$ and $y=(y_1, \ldots, y_{n+1})$ we can consider the segment
$\overline{xy}:=(1-t)x+ty$ in $\Delta$. This turns out to be a
geodesic with respect to the symmetric Lipschitz metric.  See
\cite{FrMa} for details and proofs.

However geodesics, even within a given simplex, are in general not
unique. Our strategy is broadly to determine sufficiently many
``unique'' geodesics.

\begin{dfn} A geodesic segment $\sigma$ of CV$_n$ is {\em rigid} if
for any two points on it, the restriction of $\sigma$ is the unique
(unparameterised) geodesic joining them.
\end{dfn}

We fix now a multi-theta simplex $\Delta$, and we denote by
$e_0,\dots, e_n$ the (oriented) edges of underlying graph of
$\Delta$. Any point $x$ in $\Delta$ is thus determined by lengths
$L_x(e_i)$ of $e_i$ in $x$. As usual, we denote by $\bar e_i$ the edge
$e_i$ with the inverse orientation.

We begin by describing a set of standard rigid geodesics of $\Delta$.

\begin{lemma}[Standard rigid geodesics]\label{lstandard} Let $x\neq y$
be metric graphs in $\Delta$. For any $i$ let $\lambda_i$ be the
stretching factor of $e_i$ from $x$ to $y$:
$$\lambda_i=\frac{L_y(e_i)}{L_x(e_i)}.$$
If the set of such stretching factors contains exactly two elements,
none of them with multiplicity $2$, then the segment between $x$ and
$y$ is rigid.
\end{lemma}

\proof This is a consequence of Lemma~\ref{lsausages} and
Lemma~\ref{lgeod}.

Indeed, up to rearranging edges, we can suppose
$\lambda_0=\dots=\lambda_k=\mu$ and $\lambda_{k+1}=\dots=
\lambda_n=\lambda$ for two numbers $\mu<\lambda$.  By scaling the
graph $y$ by $\mu$, we may reduce to the case where
$\lambda_0=\dots=\lambda_k=1<\lambda_{k+1}=\dots=\lambda_n=\lambda$. In
particular, we have scaled $y$ so that the edges $e_0, \ldots, e_k$
have the same length in both $x$ and $y$.

Let $z$ be a point in a geodesic joining $x$ and $y$. We claim that,
up possibly to scaling, the edges $e_0,\dots, e_k$ are not stretched
from $x$ to $z$, while the edges $e_{k+1},\dots,e_n$ are stretched all
by the same amount between $1$ and $\lambda$. That is to say, we scale
$z$ so that the length of $e_0$ in $z$ is equal to the length of $e_0$
in both $x$ and $y$. Now we claim that if $z$ belongs to a geodesic
joining $x$ and $y$, then it belongs to the segment between $x$ and
$y$, which therefore is rigid.

Let us examine our claim. First, suppose $k>1$. Then, the loops
$e_i\bar e_j$ with $i,j\leq k$ are minimally stretched from $x$ to
$y$. Thus, by Lemma~\ref{lgeod} the same must be true from $x$ to
$z$. In particular all such loops are stretched the same from $x$ to
$z$. As we have at least three such loops (because $k>1$) this implies
that the edge-stretching factors
$L_z(e_0)/L_x(e_0),\dots,L_z(e_k)/L_x(e_k)$ all coincide.

%
%
%
%

This fact is also trivially true if $k=0$, while the case $k=1$ is
impossible because the multiplicity of $\mu$ was supposed different
from $2$. So, up possibly to scaling, the edges $e_0,\dots,e_k$ are
not stretched from $x$ to $z$ (they have the same length in each
metric graph).

The same argument, now with maximally stretched loops, shows that
edges $e_{k+1},\dots,e_n$ are all stretched the same amount (as above,
$k\neq n-2$ because the multiplicity of $\lambda$ is not $2$) and by
an amount which is between $1$ and $\lambda$.  \qed

\medskip

Note that rigid segments of type just described, always emanate from
any point $x$ of $\Delta$. Indeed it suffices to consider a set $I$ of
edges and consider a point $y$ whose edge-lengths equal those of $x$
for edges in $I$ and, say, double those of $x$ for remaining edges. As
above, this will be a rigid geodesic which obviously extends to a
rigid geodesic ray. This is why we call such geodesic ``standard''.

One can think these geodesics as being a standard set in the tangent
space at $x$. Our objective now is to see that points of $\Delta$ can
have more rigid geodesics emanating from them, and that such a set of
``rigid'' directions is minimal when $x$ is the centre of $\Delta$. We
notice that this is a substantial difference with respect to case of
rose-simplices, which are homogeneous as there is transitive action of
translations.

\begin{lemma}[Rigid geodesics from the centre, in the case of rank at
least $3$]\label{lcentre} Suppose $x$ is the centre of $\Delta$. If
$n\geq3$, then any rigid geodesic through $x$ is of the type described
in Lemma~\ref{lstandard}.
\end{lemma} \proof We scale $x$ so that its edges have length one.  We
have to show that for any point $y$, if the segment $\overline{xy}$ is
rigid, then the set of edge-stretching factors contains exactly two
elements, none of them with multiplicity two.

Suppose first that we have two edge-stretching factors, one of them
with multiplicity two.  Up to scaling $y$ and rearranging edges, we
can suppose that the stretching factors of edges $e_i$ are $1$ for
$i=0,\dots,n-2$ and $\lambda$ for $i=n-1,n$. We show that in that case
the segment from $x$ to $y$ is not rigid.

Without loss of generality we can suppose $\lambda >1$.  Let $z$ be
the middle point of such segments, that is to say
$$1=L_z(e_0)=\dots=L_z(e_{n-2}) \qquad
L_z(e_{n-1})=L_z(e_n)=\frac{1+\lambda}{2}.$$

Since $\lambda>1$, the loops $e_i\bar e_j$ with $i,j<n-1$ (whose
existence is guaranteed because $n\geq3$) are minimally stretched, and
$e_{n-1}\bar e_n$ is maximally stretched, both from $x$ to $y$, form
$x$ to $z$ and from $z$ to $y$.

Moreover, since the inequalities in play are all strict, the same
remains true if we slightly perturb the length of $e_n$ (note that
maximally and minimally stretched loops have no common edges).  That
is to say, if $z_\varepsilon$ denote the graph whose edge-lengths
equal those of $z$ except for $e_n$, for which we set
$L_{z_\varepsilon}(e_n)=L_z(e_n)+\varepsilon$, for small enough
$\varepsilon$, it is still true that loops $e_i\bar e_j$ with
$i,j<n-1$ are minimally stretched, and $e_{n-1}\bar e_n$ is maximally
stretched, both from $x$ to $y$, form $x$ to $z_\varepsilon$ and from
$z_\varepsilon$ to $y$.  This
implies $$d(x,z_\varepsilon)+d(z_\varepsilon,y)=d(x,y).$$ Thus, as
segments are geodesics, the union $\sigma_\varepsilon$ of segments
$\overline{xz_\varepsilon}$ and $\overline{z_\varepsilon y}$ is a
geodesic between $x$ and $y$.  On the other hand it is clear that
$z_\varepsilon$ does not belong to segment $\overline{xy}$, so
$\sigma_\varepsilon$ is different from $\overline{xy}$ which is
therefore not rigid.

\medskip

It now remains to show that if we have at least three different
stretching factors, then we can find a geodesic between $x$ and $y$
which is not a segment.  As above, we can scale $y$, and rearrange
edges so that $1=\lambda_0\leq\lambda_1\leq\dots\leq\lambda_n$.

Since $L_x(e_i)=1$ for all $i$, the minimally stretched loops from $x$
to $y$ are all the $e_i\bar e_j$ for which $\lambda_i=\lambda_0=1$ and
$\lambda_j=\lambda_1$, and maximally stretched ones are those $e_i\bar
e_j$ for which $\lambda_i=\lambda_{n-1}$ and $\lambda_j=\lambda_n$.

Let $z$ be the middle point of the segment from $x$ to $y$.  Let
$\lambda\in\{\lambda_i\}$ be an edge-stretching factor such that
$1\neq\lambda\neq \lambda_n$. Let $z_\varepsilon$ be a metric graph
whose edge-lengths equal those of $z$, except that for edges stretched
by $\lambda$, for which differ by $\varepsilon$
$$L_{z_\varepsilon}(e_i)=
\left\{
  \begin{array}{ll} L_z(e_i) & \lambda_i\neq \lambda\\ L_z(e_i)
+\varepsilon & \lambda_i=\lambda
  \end{array} \right.
$$

and let $\sigma_\varepsilon$ be the union of segments
$\overline{xz_\varepsilon}$ and $\overline{z_\varepsilon y}$.

It is clear --- because we have at least three stretching factors ---
that $z_\varepsilon$ does not belong to the segment $\overline{xy}$,
whence $\sigma_\varepsilon\neq \overline{xy}$.  If we show that
$\sigma_\varepsilon$ is a geodesic we are done. As above, it is enough
to show that
$$d(x,z_\varepsilon)+d(z_\varepsilon,y)=d(x,y).$$

For that, we have to prove that there are loops $\gamma_0$ and
$\gamma_1$ that are respectively minimally and maximally stretched
from $x$ to $z_\varepsilon$ and from $z_\varepsilon$ to $y$. This
easily follows, for small enough $\varepsilon$, by the choice of
$\lambda$.  Indeed, it suffices (since the other cases are easier) to
look at the situation when the stretching factors are
$1,\lambda,\dots,\lambda,\lambda_n$. Here, min. and max.
lops-stretching factors form $x$ to $y$ are $(1+\lambda)/2$ and
$(\lambda+\lambda_n)/2$, realised by $e_0\bar e_i$ and $e_{i}\bar e_n$
for $i=1,\dots,n$. Such loops are therefore min and max stretched both
from $x$ to $z$ and from $z$ to $y$, and perturbing $\lambda$ a little
such loops remain min. and max. stretched.  \qed

\medskip

Now, we show how Lemma~\ref{lcentre} provides (in rank bigger than
two) a metric characterisation of the centre of $\Delta$ as the point
having the minimum number of rigid geodesics passing through it.

\begin{lemma}\label{lnoncentre} For any point $x$ other than centre of
$\Delta$, there is at least one rigid geodesic emanating from $x$
which is not of the type described in Lemma~\ref{lstandard}.
\end{lemma}

\proof We denote by $x_i$ the lengths $L_x(e_i)$. Up to scaling $x$
and rearranging edges, we can suppose that
$$1=x_0\geq x_1\geq\dots\geq x_n.$$ We want to find stretching factors
$$1=\lambda_0\leq \lambda_1\leq\dots\leq \lambda_n$$ at least three of them
being different, such that segment between $x$ and point $y$
corresponding to graph whose edges have length $\lambda_i x_i$, is
rigid. As three of the $\lambda_i$ are different, this will prove the
lemma.

Let us start by making the simplifying assumption that $x_n\neq x_1$.

Stretching factors, from $x$ to $y$, of loops $e_i\bar e_j$ are
$\frac{\lambda_i x_i+\lambda_j x_j}{x_i+x_j}$, and if $1=\lambda_0\leq
\lambda_1\leq\dots\leq\lambda_n$, an immediate calculation shows that
whenever $j\geq i$ we have
$$\frac{1+\lambda_i x_i}{1+x_i}\leq
\frac{\lambda_i x_i+\lambda_j x_j}{x_i+x_j} \leq \frac{\lambda_j
x_j+\lambda_n x_n}{x_j+x_n}.
$$

 This implies that if we are searching for minimally (respectively
maximally) stretched loops, we can restrict to loops of the form
$e_0\bar e_i$ (respectively $e_i\bar e_n$.)

The idea is now to force such loops to have the same stretching
factors. We impose conditions
$$ \lambda_1=\frac{1+x_1}{x_1}$$
and, for $i>0$

\begin{equation}
  \label{esf} \frac{1+\lambda_i x_i}{1+x_i} = \frac{1+\lambda_1
x_1}{1+x_1} =\frac{2+x_1}{1+x_1}
\end{equation}

We remark that the assumption on $\lambda_1$ is for simplifying
calculations, we only need $\lambda_1>1$.

We can solve these equations getting
$$\lambda_i=
\frac{(2+x_1)(1+x_i)}{(1+x_1)x_i} -\frac{1}{x_i}=
1+\frac{1+x_i}{(1+x_1)x_i}
$$
thus $\lambda_i\geq\lambda_1$, with equality if and only if $x_i=x_1$,
and $\lambda_i\leq\lambda_j$ for $j\geq i$, with equality if and only
if $x_i=x_j$. In particular, under our simplifying assumption, we have
$\lambda_0=1<\lambda_1<\lambda_n$, so at least three of the
$\lambda_i$'s are different.

So we get numbers $\lambda_i$'s with the requested properties.  Now,
let $y$ be the point of $\Delta$ given by
$$L_y(e_i)=\lambda_i x_i$$

and let $z$ be any point in a geodesic between $x$ and $y$, scaled so
that $L_z(e_0)=1$. We define $\mu_i$ by
$$L_z(e_i)=\mu_i x_i.$$

Loops $e_0\bar e_i$ are minimally stretched from $x$ to $y$.  Thus, we
must have that such loops are minimally stretched from $x$ to $z$ and
from $z$ to $y$. This forces the edge-stretching factors $\mu_i$ to
satisfy condition $(\ref{esf})$, which allows us to obtain $\mu_i$ as
a function of $\mu_1$ exactly as $\lambda_i$ is obtained from
$\lambda_1$. This implies that, if $z'$ is the point in the geodesic
line between $x$ and $y$ with first edge-stretching factor equal to
$\mu_1$, we have that $z=z'$.

So that $z$ belongs the segment $\overline{xy}$ which is hence rigid,
and not of the type described in Lemma~\ref{lstandard}.

We are now left with the case in which $x_n=x_1$ so $x_i=x_j$ for any
$i,j\neq 0$. As we are supposing that $x$ is not the centre of
$\Delta$, we must have $x_0\neq x_1$. Up to scaling $x$ and
rearranging edges, this case is equivalent to
$$(x_0,\dots,x_n)=(1,\dots,1,c)$$
with $c>1$.

We choose $y$ of the form
$$y=(1,\lambda,\dots,\lambda,\mu c)$$
Stretching factors of simple loops are
$$\frac{1+\lambda}{2},\qquad \frac{1+\mu c}{1+c},\qquad
\lambda,\qquad \frac{\lambda+\mu c}{1+c}.$$ Now, we impose conditions
$$\mu c=1, \qquad \frac{1+\lambda}{2}=\frac{1+\mu c}{1+c}$$
which imply that $\lambda\neq\mu$ because $c\neq 1$, and
$\lambda<1$. Whence
$$\frac{1+\lambda}{2}=\frac{1+\mu
  c}{1+c}>\max(\lambda,\frac{\lambda+\mu c}{1+c}).$$

So all the loops $e_0\bar e_i$ are maximally stretched from $x$ to $y$
(and in particular, stretched by the same amount).  Now we argue as
before: the same must be true for any point $z$ on any geodesic from
$x$ to $y$, and this forces $z$ to be of the form (once scaled so that
$L_z(e_0)=1$)
$$z=(1,\bar\lambda,\dots,\bar\lambda,\bar\mu c)$$
with
$$\frac{1+\bar \lambda}{2}=\frac{1+\bar \mu c}{1+c}.$$

As above, this implies that $z$ belongs to the segment
$\overline{xy}$, which is then rigid and it is not of the type
described in Lemma~\ref{lstandard} because
$1\neq\lambda\neq\mu\neq1$.\qed

\begin{lemma}[Rigid geodesics in rank two]\label{lrigeod2} Let $x\neq
y$ be two marked metric graphs in $\Delta$. Suppose $n=2$, so that
$\Delta$ has exactly three different (unoriented) simple loops.  Then
the segment $\overline{xy}$ is rigid if and only if two of the three
simple loops are stretched the same from $x$ to $y$.
\end{lemma} \proof The proof use same arguments of higher rank case,
but takes in account the peculiarities of rank two.

If the three simple loops are stretched by three different factors,
then for any point $w$ close enough to the middle point $z$ of
$\overline{xy}$, the maximally and minimally stretched loops do not
change from $x$ to $w$ from $w$ to $y$ and from $x$ to $y$. So that
$\overline{xy}$ is not rigid.

On the other hand, if two simple loops are stretched by the same
factor, we may rearrange the edges so that $e_0$ is the edge shared by
such loops, and scale graphs so that $L_x(e_0)=L_y(e_0)=1$. Moreover,
as we have only three simple loops, $e_0\bar e_1$ and $e_0\bar e_2$
are either maximally or minimally stretched from $x$ to $y$. So the
same must be true from $x$ to $z$ and from $z$ to $y$ for any point
$z$ in a geodesic between $x$ and $y$.  If $x=(1,a,b)$ and
$y=(1,\lambda a,\mu b)$, we have
$$\frac{1+\lambda a}{1+a}=\frac{1+\mu b}{1+b}$$
and the same relation holds for the edge stretching factors of point
$z$ which therefore belongs to the segment $\overline{xy}$.\qed

\begin{lemma}\label{lucr} For any rose-face of $\Delta$ there is a
unique rigid geodesic from the centre of $\Delta$ to that face.
\end{lemma} \proof By Lemma~\ref{lcentre} and~\ref{lrigeod2}, a rigid
geodesic emanating from the centre is of the type described in
Lemma~\ref{lstandard} (and \ref{lrigeod2} in the rank-$2$ case). A
rose-face corresponds to collapsing an edge, say $e_0$. So in a rigid
geodesic from the centre to that face we have $\lambda_0=t$ and
$\lambda_i=1$ for $i>0$, with $t\in[1,0]$. Therefore such geodesic is
unique.\qed

\medskip

Now, we continue with proof of Theorem~\ref{tmultitheta}. We begin by
examining the first claim in the rank-two case.  Since permutations of
edges of $\Delta$ are isometries that fix its centre and permute its
rose-faces, up to composing $\Phi$ with such a permutation we can
suppose that $\Phi$ does not permute rose-faces of $\Delta$. If the
restriction of $\Phi$ to a rose-face has a translational part, then
for any point $x$ in that face we see that the distance of $\Phi^n(x)$
from at least one of the remaining two rose faces of $\Delta$ goes to
infinity, this being impossible because $\Phi$ is an isometry. It
follows that $\Phi$ fixes the centres of rose-faces of
$\Delta$. Explicit calculations (using Lemma~\ref{lrigeod2}, see the
Appendix) show that the centre of $\Delta$ is the unique point which
is joined to the centres of the three rose-faces by rigid
geodesics. Thus $\Phi$ fixes the centre of $\Delta$, and the first
claim of Theorem~\ref{tmultitheta} is proved for $n=2$.

If $n\geq 3$, Lemma~\ref{lcentre} and Lemma~\ref{lnoncentre} imply
that any isometry $\Phi$ of $\Delta$ must fix its centre, so first
claim of theorem is proved. Moreover, if $\Phi$ does not permute
rose-faces, then by Lemma~\ref{lucr} it must fix point-wise rigid
geodesics emanating from $x$ and going to rose-faces. In particular,
$\Phi$ fixes centres of rose-faces.

\begin{rem}\label{rrn} Note that we have proved that if $\Phi$ does
not permute rose-faces of $\Delta$, then its restriction to any
rose-face has no translational parts, which is to say that it fixes
the centre of rose-face.
\end{rem}

Therefore, by Theorem~\ref{tisorose}, restriction of $\Phi$ to
rose-faces of $\Delta$ is an element of $
S_n\times\langle\sigma\rangle$. In the next lemma we show that such an
element must be the identity. We first introduce some terminology.

Let $R_i$ denote the rose-face of $\Delta$ obtained by collapsing edge
$e_i$, and let $C_i$ denote its centre.  Also, for $i\neq j$ we let
$\Gamma_j^i(\epsilon)$ denote the point of $R_j$ all of whose petals
have length $1$ except for $e_i$ which has length $\epsilon$. For
$i\neq j$, straightforward calculations show

\begin{equation}\label{e6} d(C_i, \Gamma_j^k(2)) \ = \ \left\{
\begin{array}{lclcl} \log 6 &~& i=k &~& \textrm{in any rank}\\ \log 6
&~& i\neq k &~& \textrm{in rank bigger that } 2\\ \log 3 &~& i\neq k
&~& \textrm{in rank } 2

\end{array} \right.
\end{equation} and

\begin{equation}\label{e7} d(C_i, \Gamma_j^k(0.5)) \ = \ \left\{
\begin{array}{lclcl} \log 3 &~& i=k &~& \textrm{in any rank}\\ \log 8
&~& i\neq k &~& \textrm{in rank bigger that } 2\\ \log 6 &~& i\neq k
&~& \textrm{in rank } 2
\end{array} \right.
\end{equation}

Note that in the rank-$2$ case, for $i\neq j\neq k$ we have
$\Gamma_j^k(2)=\Gamma_j^i(0.5)$, up to scaling.
\begin{lemma}\label{lsigma} Let $\Phi$ be an isometry of a multi-theta
simplex $\Delta$ which fixes the centres of its rose-faces. Then
$\Phi$ is the identity on each rose face.
\end{lemma} \proof By Theorem~\ref{tisorose}, the restriction of
$\Phi$ to any rose face is an element of $ S_{n}\times\langle
\sigma\rangle$.  Hence the image of the point $\Gamma_j^i(2)$ is
either $\Gamma_j^{k}(2)$ or $\Gamma_j^{k}(0.5)$ for some $k$.
However, by $(\ref{e6})$ and $(\ref{e7})$, since each $C_i$ is fixed
by $\Phi$ the distances to $\Gamma_j^i(2)$ are preserved and we must
have that $\Gamma_j^i(2)$ is actually fixed by $\Phi$.  Since this is
true for every $i\neq j$, and since the only element of
$S_{n}\times\langle\sigma\rangle$ which fixes all these is the
identity, we get that $\Phi$ restricts to the identity on any $R_j$.
\qed

\medskip

We can now finish proof of Theorem~\ref{tmultitheta}.  We proved that
any isometry of $\Delta$ fixes it centre, and that if it does not
permute rose-faces $R_i$, then it fixes their centres $C_i$. By
Lemma~\ref{lsigma} this implies that $\Phi$ point-wise fixes
rose-faces of $\Delta$. Now let $x\in\Delta$. For any $y$ in some
$R_i$, we have $d(x,y)=d(\Phi(x),y)$ (because $R_i$ is fixed.)
Therefore, Theorem~\ref{tasymptotic} tells us that lengths of
a.s.c.c. in $x$ and $\Phi(x)$ coincide. Thus, by Lemma~\ref{lsausages}
we have $d(x,\Phi(x))=0$.  It follows that $\Phi$ is the identity of
$\Delta$, and the proof of Theorem~\ref{tmultitheta} is concluded.\qed

\bigskip

We come now the other main result of this section, that is that for an
isometry of CV$_n$, what happens on a single multi-theta simplex
determines the isometry on the whole CV$_n$.  The first step is to
show that if an isometry of a multi-theta simplex is the identity on a
rose-face, then it is the identity of the multi-theta simplex. Our
claim will follow then by an argument of connection.

\begin{lemma}
\label{extend} Let $\Phi$ be an isometry of a multi-theta simplex
$\Delta$ which restricts to the identity on one of the rose-face of
$\Delta$.  Then $\Phi$ restricts to the identity on each rose-face of
$\Delta$.
\end{lemma} \proof Let $R_0$ be the rose-face fixed by hypothesis.  By
Lemma~\ref{lsigma}, it is sufficient to show that $\Phi$ fixes each
centre $C_i$. By first claim of Theorem~\ref{tmultitheta}, we know
that the centre of $\Delta$ is fixed. By Lemma~\ref{lucr}, there is a
unique rigid geodesic from the centre to each rose face, ending in
$C_i$. Hence, the $C_i$ are permuted by $\Phi$.

However, the stabiliser in $Out(F_n)$ of $\Delta$ contains a subgroup
isomorphic to $S_{n+1}$, by simply permuting the edges of the
underlying graph of $\Delta$, and this subgroup will induce every
permutation of the $n+1$ rose-faces of $\Delta$. Hence, by
Lemma~\ref{lsigma}, $\Phi$ is equal to the restriction of some element
of $S_{n+1}$ (in fact, some such element which fixes the edge
corresponding to the fixed rose-face). But the only element of this
sort which restricts to the identity in a rose face is the identity.
(This also follows from $(\ref{e6})$ and $(\ref{e7})$).

\qed

\begin{thm}\label{troses} Let $\Phi$ be an isometry of CV$_n$ that
point-wise fixes a multi-theta simplex. Then it point-wise fixes all
rose and multi-theta simplices of CV$_n$.
\end{thm} \proof

We start by doing a simple calculation. Let $\Delta$ be a multi-theta
simplex of CV$_n$, with edges oriented and labelled $e_0, e_1, \ldots,
e_n$. For an edge $e$, we denote by $\overline e$ the edge $e$ with
inverse orientation. Let $R_i$ be rose face of $\Delta$ obtained by
collapsing $e_i$.  We will label the edges of $R_i$, $e_0^i, e_1^i,
\ldots, e_{i-1}^i, e_{i+1}^i, \ldots, e_n^i$.

Now let us explicitly write down the homotopy equivalences between
$R_0$ and $R_i$ in terms of these edges. The map from $R_0$ to $R_i$
is given by the following,
\begin{equation}
\label{he1}
\begin{array}{rcl} e_j^0 & \mapsto & e_j^i \overline{e_0^i}, j \neq i
\\ e_i^0 & \mapsto & \overline{e_0^i}.
\end{array}
\end{equation}

Similarly, the map from $R_i$ to $R_0$ is given by,
\begin{equation}
\begin{array}{rcl} e_j^i & \mapsto & e_j^0 \overline{e_i^0}, j \neq 0
\\ e_0^i & \mapsto & \overline{e_i^0}.
\end{array}
\end{equation}

This in particular implies that the sub-complex of CV$_n$ consisting
of multi-theta and rose-simplices is connected, as we realised Nielsen
automorphisms passing from a rose-face to another in a multi-theta
simplices.

Now, it would seem that we are done simply by starting from our
initial fixed multi-theta simplex and extending our results, via
Lemma~\ref{extend}, over the whole of CV$_n$. The only problem is that
we do not know, {\em a priori}, that $\Phi$ does not induce some
non-trivial permutation of the multi-theta simplices. Therefore, we
need to rule out this possibility.

\begin{rem} The next Lemma is an ``elementary'' proof of the fact that
permutations of multi-theta simplices do not occur. The calculations
it involves are somewhat tedious and the reader may prefer to invoke
the result of Bridson and Vogtmann~\cite{BrVo01} asserting that
simplicial actions on the spine of CV$_n$ (see~\cite{BrVo01} for
definitions and details) come from automorphisms (for $n\geq 3$).
Then, she could show that isometries naturally induce such actions on
the spine, and since the spine encodes the combinatoric of roses and
multi-theta incidences, get the desired result. We present here the
proof of Lemma~\ref{noperm} as follows because it is self-contained
and more in the spirit of the techniques of the present work.
\end{rem}

\begin{lemma}
\label{noperm} Let $\Delta$ be a multi-theta simplex, $R$ a rose face
of it. Suppose that $\Delta_1, \ldots, \Delta_k$ are all the other
multi-theta simplices in CV$_n$ which are incident to $R$. Let $\Phi$
be an isometry of CV$_n$ which point-wise fixes $\Delta$ (and
therefore $R$). Then $\Phi$ leaves each $\Delta_i$ invariant.
\end{lemma}

\proof Consider our multi-theta simplex $\Delta$ which is given by a
graph with $2$ vertices and $n+1$ edges, ordered and labelled $e_0,
\ldots, e_{n+1}$. As usual, for an edge $e$ we denote by $\overline e$
the one with inverse orientation. Moreover, we chose orientations so
that the $e_i$'s share the same initial vertex (so they also share the
terminal vertex). We will let $R$ denote the rose simplex obtained by
collapsing the edge $e_0$.

It is now an easy exercise to see that there are $2^{n-1}$ multi-theta
simplices incident to $R$. Therefore, the result is trivial in CV$_2$
and we shall restrict our attention to CV$_n$ for $n \geq 3$.

We shall describe the set of multi-thetas incident to $R$ by listing
the homotopy equivalences from $\Delta$. Specifically, choose some $I
\subseteq \{ 1, \ldots, n \}$ and consider the homotopy equivalence on
$\Delta$ given by,
$$
\begin{array}{rcl} e_0 & \mapsto & e_0 \\ e_i & \mapsto & e_i, \ i
\not\in I \\ e_i & \mapsto & e_0 \overline{{e_i}} e_0, \ i \in I\\
\end{array}
$$

It is then clear that the set of all multi-thetas incident to $R$ will
be given by these maps. However, we note that replacing $I$ by its
complement gives the same simplex, so we have counted each twice. From
now we will make a choice between $I$ and $I^c$ so that $| I | \geq |
I^c |$ (or $I=\emptyset$) --- if $|I |= | I^c |$ the choice will be
arbitrary. Hence if $I$ is not empty it will have at least two
elements, and its complement will be non-empty. Let $\Delta_I$ denote
the multi-theta simplex obtained via the map above. This gives us our
$2^{n-1}$ multi-thetas, with $\Delta=\Delta_{\emptyset}$.

Now, we will show that the distances from $\Delta$ will determine the
$\Delta_I$. Note that since we are dealing with multi-theta graphs, by
the Sausages Lemma (\ref{lsausages}) the maximally and minimally
stretched loops can be taken to be simple closed curves, which are
straightforward to enumerate. Below, we present a list of curves. On
the left side, we have curves in $\Delta$ and on the right side their
image in $\Delta_I$ so that each simple closed curve in either
$\Delta$ or $\Delta_I$ appears somewhere on the list (up to
orientation). Throughout, we have that $i, j \neq 0$.
$$
\begin{array}{rcl} \Delta & & \Delta_I \\ e_i \overline{{e_0}} &
\mapsto & e_i \overline{{e_0}}, \ i \not\in I \\ e_i \overline{{e_0}}
& \mapsto & e_0 \overline{{e_i}}, \ i \in I \\ e_i \overline{{e_j}} &
\mapsto & e_i \overline{{e_j}}, \ i, j \not\in I \\ & \mapsto & e_0
\overline{e_i} {{e_j}} \overline{{e_0}}, \ i,j \in I \\ & \mapsto &
e_0 \overline{ {e_i} } e_0 \overline{{e_j}}, \ i \in I, j \not\in I \\
\overline{ {e_i} } e_0 \overline{{e_j}} e_0 & \mapsto & {e_i}
\overline{{e_j}}, \ i \in I, j \not\in I \\
\end{array}
$$

Now let us assign edge lengths and calculate distances. For each
$\Delta_I \neq \Delta$, we will let all edge lengths equal $1$, since
we know that isometries preserve the centres. Next choose some $J
\subseteq \{ 1, \ldots, n \}$ and let $\Delta(J, 1/3)$ be the graph
$\Delta$ where each edge has length 1 except for the $e_j$ which has
length $1/3$ for all $j \in J$. Moreover, let us stipulate that $J
\neq \emptyset, \{ 1, \ldots, n \}$. It is then an easy exercise to
check the stretching factors for each of the simple loops in
$\Delta(J, 1/3)$ and $\Delta_I$. Clearly, this depends on the
relationship between $J$ and $I$. We list below, the possible
stretching factors between $\Delta(J,1/3)$ and $\Delta_I$, with the
condition which allows it. Some stretching factors can occur in more
than one way, in which case we have removed the redundancy (an empty
condition means the stretching factor is always realisable).

So distance is computed by taking the log of the ratio of the maximum
over the minimum of the allowed factors. Recall that by the choice we
made for $I$, we always have $|I|\geq2$ and $|I^c|\geq 1$.

$$
\begin{array}{rl} \mbox{\rm Stretching factor} & \mbox{\rm Condition}
\\ 1 & \\ 3/2 & \\ 3 & | I \cap J | \geq 2 \ \mbox{\rm or } \ |I^c
\cap J | \geq 2 \\ 2 & | I \cap J^c | \geq 1 \ \mbox{\rm and } \ |I^c
\cap J^c | \geq 1 \\ 3 & | I \cap J | \geq 1 \ \mbox{\rm and } \ |I^c
\cap J^c | \geq 1 \\ 3 & | I \cap J^c | \geq 1 \ \mbox{\rm and } \
|I^c \cap J | \geq 1 \\ 6 & | I \cap J | \geq 1 \ \mbox{\rm and } \
|I^c \cap J | \geq 1 \\ 1/2 & | I \cap J^c | \geq 1 \ \mbox{\rm and }
\ |I^c \cap J^c | \geq 1 \\ 3/5 & | I \cap J | \geq 1 \ \mbox{\rm and
} \ |I^c \cap J^c | \geq 1 \\ 3/5 & | I \cap J^c | \geq 1 \ \mbox{\rm
and } \ |I^c \cap J | \geq 1 \\ 3/4 & | I \cap J | \geq 1 \ \mbox{\rm
and } \ |I^c \cap J | \geq 1 \\
\end{array}
$$

We now apply these conditions to calculate the distances from
$\Delta(J, 1/3)$ to $\Delta_I$ when $J$ has exactly 2
elements. Specifically,
\begin{itemize}
\item If $|J|=2$, $|J \cap I| = |J \cap I^c|=1$, and $|I^c| \geq 2$,
then maximal and minimal stretching factors are $6$ and $1/2$, so the
distance is $\log 12$.

\item If $|J|=2$, $|J \cap I| = |J \cap I^c|=1$, and $|I^c|=1$, then
the max and min stretching factors are $6$ and $3/5$, whence the
distance is $\log 10$.

\item If $|J|=2$ and $J \subseteq I$ or $J \subseteq I^c$, then the
maximal stretching factor is always $3$, and the distance is $\log5$
or $\log 6$, depending on the sizes of $I$ and $I^c$.
\end{itemize}

Hence, we may determine $I$ and $I^c$. More precisely, the set
$\{1\}\cup\{i\neq 1:d(\Delta(\{1,i\}, 1/3),\Delta_I)=\log 5\} \cup
\{i\neq 1:d(\Delta(\{1,i\}, 1/3),\Delta_I)=\log 6\}$ is equal to
either $I$ or $I^c$.

Note that this doesn't let us distinguish which one we picked, but
since $\Delta_I$ only depended on the pair $I, I^c$, this is
sufficient to distinguish the simplex and proves Lemma~\ref{noperm}.
\qed

\medskip

Now Theorem~\ref{troses} follows.\qed

\section{Proof of Theorem~\ref{tmain} and main results}

We prove here results stated in Section~\ref{sintro}.

\proof[Proof of Theorem~\ref{tmain}] Our claim is that the
isometry-group of outer space of rank-$n$ free group, is just
$Out(F_n)$ (for $n \geq 3$). Clearly, $Out(F_n)$ acts faithfully on
CV$_n$ for $n \geq 3$ and this action is by isometries (see for
instance~\cite{FrMa}). Thus, we have an inclusion of $Out(F_n)$ into
the group of isometries of CV$_n$. For $n=2$, we still have an
isometric action, but this is no longer faithful. However, up to this
small kernel (a group of order $2$ consisting of the identity and the
automorphism which inverts each basis element), we still have a map
from $Out(F_2)$ to the isometry group of CV$_2$.

Our goal is to show that this exhausts the isometry group of CV$_n$
(in either case).

Let $\Phi$ be an isometry of CV$_n$. We shall compose $\Phi$ with
elements of $Out(F_n)$ until we obtain the identity.

By Lemma~\ref{L2}, $\Phi$ maps multi-theta simplices to multi-theta
simplices. Therefore, since the action of $Out(F_n)$ on multi-theta
simplices is transitive, we may suppose that $\Phi$ leaves invariant a
multi-theta simplex $\Delta$. In fact, the stabiliser in $Out(F_n)$ of
$\Delta$ will induce any permutation of the $n+1$ rose faces of
$\Delta$ and so we may also assume that $\Phi$ leaves both $\Delta$
and every rose-face of $\Delta$ invariant.

Theorem~\ref{tmultitheta} then implies that $\Phi$ is the identity of
$\Delta$. Then, by Theorem~\ref{troses} $\Phi$ point-wise fixes all
rose-simplices. And Theorem~\ref{tid} implies that $\Phi$ is the
identity.\qed

\medskip

\proof[Proof of Theorem~\ref{tm2}] Let us denote by
$\textrm{Isom}_R(CV_n)$ the group of isometries of CV$_n$ for the
non-symmetric metric $d_R$, and by $\textrm{Isom}(CV_n)$ the group of
isometries of CV$_n$ for the symmetric metric $d$.

Let $\Phi$ be an isometry of CV$_n$ for $d_R$. Then
$$\forall x,y\qquad d_R(x,y)=d_R(\Phi x,\Phi y).$$
Since $d_L(x,y)=d_R(y,x)$, we have that $\Phi$ is also an isometry for
$d_L$, whence $\Phi$ is an isometry for the symmetric Lipschitz metric
$d$.  Thus, $\textrm{Isom}_R(CV_n)\subseteq\textrm{Isom}(CV_n)$.

As for the symmetric case, one has that $Out(F_n) \subseteq
\textrm{Isom}_R(CV_n)$ (with a small adjustment for rank $2$).  By
Theorem~\ref{tmain} we have that Isom$(CV_n)=Out(F_n)$. Thus
$$Out(F_n) \subseteq Isom_R(CV_n)\subseteq
Isom(CV_n)=Out(F_n).$$ The same for $d_L$. (The argument for rank $2$
is the same.)  \qed

\medskip
\proof[Proof of Corollary~\ref{c}]
Every homomorphism from $\Gamma$ to $Out(F_n)$ has finite image
by~\cite{BrFarb}, and every finite subgroup of $Out(F_n)$ has a fixed
point in its action on $CV_n$ by~\cite{CU}.
\qed

\appendix
\section{Rigid geodesics in rank two}

Here we explicitly calculate rigid geodesics emanating from centres of
rose-simplices and pointing into theta-simplices, for the rank-two
case. Showing that for any theta-simplex, its centre is the unique
point simultaneously joined to centres of all rose-faces by rigid
geodesics.

We fix a theta-simplex $\Delta$ and we parametrise its points by
(projective classes of) triples of positive numbers $(x,y,z)$. Such
simplex is a triangle with vertices removed, as can be seen by taking
representatives unitary volume.

Let $(1,0,1)$ be the centre of a rose face of $\Delta$ and let
$(1,z,y)$ be a point joined to it by a rigid segment, scaled so that
$x=1$. Stretching factors are
$$1+z,z+y,\frac{1+y}{2}$$
(the loop with stretching factor $z+(1+y)/2$ is not relevant)

By Lemma~\ref{lrigeod2} we must have only two stretching factors from
$(1,0,1)$ to $(1,z,y)$. Possible cases are $1+z=z+y,
1+z=\frac{1+y}{2}, z+y=\frac{1+y}{2}$. If

$$1+z=z+y$$
then $y=1$ and this is the rigid geodesic going to the centre of
$\Delta$. If

$$1+z=\frac{1+y}{2}$$
then $y=1+2z$, then $(1,y,z)=(1,z,1+2z)$. We want to know where such
geodesic hits other rose-faces.  Letting $z\to\infty$ and scaling by
$z$ we get $(1/z,1,2+1/z)$ which ends up to the point $(0,1,2)$.
Finally,
$$z+y=\frac{1+y}{2}$$
gives $z=(1-y)/2$, so that $(1,y,z)=(1,(1-y)/2,y)$. Letting $y\to 0$
we get $(1,0.5,0)$. The picture of rigid geodesics through the centres
is therefore as follows

\setlength{\unitlength}{1ex}
  \begin{picture}(65,35) \put(15,5){\line(3,5){15}}
\put(15,5){\line(1,0){30}} \put(45,5){\line(-3,5){15}}
\put(30,5){\line(0,1){10}} \put(22.5,17.5){\line(3,-1){7.5}}
\put(22.5,17.5){\line(3,1){12.5}} \put(22.5,17.5){\line(1,-3){4.15}}
\put(37.5,17.5){\line(-3,-1){7.5}} \put(37.5,17.5){\line(-3,1){12.5}}
\put(37.5,17.5){\line(-1,-3){4.15}} \put(30,5){\line(1,1){9.35}}
\put(30,5){\line(-1,1){9.35}}

    \put(15,17.5){$(1,0,1)$} \put(35,22){$(0,1,2)$}
\put(20,2){$(1,0.5,0)$}

    \put(10,30){$z=0$} \put(17,30){\vector(2,-1){9}}

    \put(44,30){$x=0$} \put(43,30){\vector(-2,-1){9}}

    \put(38,0){$y=0$} \put(41,2){\vector(0,1){2.5}}
  \end{picture}

\end{document}